\documentclass{article}
\usepackage[english]{babel}
\usepackage{amssymb,amsmath}
\usepackage{upgreek}
\usepackage{graphicx,float}
\usepackage{url}
\newcommand{\R}{\mathbb{R}}

\newcommand{\Z}{\mathbb{Z}}
\newcommand{\WW}{\textbf{W}}
\newcommand{\FF}{\textbf{F}}
\newcommand{\TT}{\textbf{T}}
\newcommand{\NN}{\textbf{N}}
\newcommand{\BB}{\textbf{B}}

\newcommand{\U}{{\textbf{U}}}
\newcommand{\UM}{{\textbf{UM}}}
\newcommand{\UP}{{\textbf{UP}}}
\newcommand{\ii}{{\textbf i}}
\newcommand{\ji}{{\textbf j}}
\newcommand{\ki}{{\textbf k}}
\newcommand{\nn}{{\textbf n}}
\newcommand{\mm}{{\textbf m}}
\newcommand{\VV}{\textbf{V}}
\newcommand{\SE}{\mathcal{\textbf{S}}}
\newcommand{\fig}{F{\sc{ig}}. }

\newcommand{\vectrois}[3]{\left(\begin{array}{c}#1\\#2\\#3\end{array}\right)}
\newtheorem{thm}{Theorem}[section]
\newtheorem{defi}{Definition}[section]
\newtheorem{lem}{Lemma}[section]
\newtheorem{rque}{Remark}[section]

\def\abs#1{\vert #1 \vert}
\def\At{\widetilde{A}}
\def\ct{c(\widetilde \WW)}
\def\ctt{\widetilde c}
\def\jump#1{\left[\left[ #1 \right]\right]}
\def\dive{{\rm div}}
\def\dsp{\displaystyle}

\def\Qt{\widetilde{Q}}
\def\ut{\widetilde{u}}

\begin{document}
\title{A model for unsteady mixed flows in non uniform closed water pipes and a well-balanced finite volume scheme}
\author{Christian Bourdarias\thanks{email: Christian.Bourdarias@univ-savoie.fr,},
Mehmet Ersoy\thanks{email: Mehmet.Ersoy@univ-savoie.fr,}~and
St\'{e}phane Gerbi\thanks{email: Stephane.Gerbi@univ-savoie.fr}\\[0.3cm]
Laboratoire de Math\'ematiques, Universit\'e de Savoie\\
73376 Le Bourget du Lac, France}
\date{}

\maketitle
\begin{abstract}
We present the derivation of a new unidirectional model for unsteady mixed flows in non uniform closed water pipes. We introduce a local reference frame  to take into account the local perturbation caused by the changes of section and slope. Then an asymptotic analysis is performed to obtain a model for free surface  flows and another one for  pressurized flows. By coupling these  models through the transition points  by the use of a common set of variables and a suitable pressure law, we obtain a simple formulation   called PFS-model close to the shallow water equations with source terms. It takes into account  the changes of section and the slope variation in a continuous way through transition points. Transition point between the two types of flows is treated as a free boundary associated to a discontinuity of the gradient of pressure.  The numerical simulation is  performed by making use of a Roe-like finite volume scheme that we adapted to take into account geometrical source terms in the convection matrix. Finally some numerical tests are presented.
\end{abstract}
\section{Introduction}
The presented work takes place in a more general framework: the modelling of unsteady mixed flows
in any kind of closed pipe taking into account  the cavitation problem and air entrapment. We are
interested in flows occurring in closed pipes with non uniform sections, where some parts of the flow
can be free surface (it means that only a part of the pipe is filled) and other parts are
pressurized (it means that the pipe is full). The transition phenomenon between the two types of
flows occurs in many situations such as storm sewers, waste or supply pipes in hydroelectric
installations. It can be induced by sudden changes in the boundary conditions as failure pumping.
During this process, the pressure can  reach severe values and  may cause damages.
The simulation of such a phenomenon is thus a major challenge and a great amount of works was devoted to it these last
years (see \cite{CSZ97},\cite{N90},\cite{R85},\cite{SWB98}, for instance).

The classical shallow water equations are commonly used to describe  free surface flows in open channels. They are also used in the study of mixed flows  using the Preissman slot artefact (see for example \cite{CSZ97,SWB98}). However, this technic does not take into account  the  depressurisation phenomenon which occurs during a water hammer except in recent works \cite{KAEDP09,KAEDM08,KDAEDM08}. On the other hand  the Allievi equations, commonly used to describe pressurized flows, are written in a non-conservative form which is not well adapted to a natural coupling with the shallow water equations.

A model for the unsteady mixed water flows in closed pipes and a finite volume
discretisation have been previously studied by two of the authors \cite{BG07} and a kinetic formulation  has
been proposed in \cite{BGG08}. We propose here the \textbf{PFS}-model which tends to extend naturally the work in \cite{BG07} in the case of a closed pipe with  non uniform section. For the sake of simplicity, we do not deal with the deformation of the domain induced by the change of pressure. We will consider only an infinitely rigid pipe.

The paper is organized as follows.
Section 2 is devoted to the derivation of  the free surface model from the $3$D \emph{incompressible} Euler equations which are written in a
suitable local reference frame (following \cite{BFML07,BMCPV03}) in order to take into account the local effects produced by the changes of
section and the slope variation. The construction of the free surface model is done by a formal asymptotic analysis.
Seeking for an approximation at first order gives the model called \textbf{FS}-model.
In Section~\ref{SectionPFModel}, we adapt the derivation of the \textbf{FS}-model to derive the pressurized model, called \textbf{P}-model, from the $3$D \emph{compressible} Euler equations.
Writing the source terms of these two models, \textbf{P} and \textbf{FS}-model, into a unified form and  using the same couple of conservative unknowns as in \cite{BG08}, we propose  in Section \ref{SectionCoupling} a model  for mixed flows, that we call \textbf{PFS}-model . We state some mathematical properties of this model. Section \ref{PFSRoe} is
devoted to the extension of the VFRoe  scheme described in  \cite{BGH00,GHS01,BG07} that was used for the case of uniform pipes. In Section \ref{SectionWB}, we show how to construct a  convection matrix in order to get an exactly well-balanced scheme. Several numerical tests are presented in Section \ref{NumericalTests}.

\newpage
{\scriptsize
\subsubsection*{Notations concerning geometrical variables}

\begin{itemize}
\item $(0,\ii,\ji,\ki)$: cartesian reference frame 
\item $\omega(x,0,b(x))$: parametrization in the reference frame $(0,\ii,\ji,\ki)$ of the plane curve $\mathcal{C}$ which corresponds to the main flow axis
\item $(\TT,\NN,\BB)$: Serret-Frenet reference frame attached to $\mathcal{C}$ with $\TT$ the tangent vector, $\NN$ the normal vector
and $\BB$ the binormal vector
\item $X,Y,Z$: local variable in the Serret Frenet reference frame with $X$ the curvilinear abscissa, $Y$ the width of pipe, $Z$ the
$\BB$-coordinate of any particle.
\item $\sigma(X,Z)=\beta(X,Z)-\alpha(X,Z)$: width of the pipe at altitude $Z$ with $\beta(X,Z)$ (resp. $\alpha(X,Z)$) is
 the Y-coordinate of right (resp. left) boundary point at altitude $Z$
\item $\theta(X)$:  angle $(\ii,\TT)$
\item $S(X)$:  cross-section area
\item $R(X)$: radius of the cross-section $S(X)$
\item $\nn_{\textbf{wb}}$: outward normal vector to the wet part of the pipe
\item $\nn$: outward normal vector at the boundary point $m$ in the  $\Omega$-plane defined below
\end{itemize}
\subsubsection*{Notations concerning the free surface (FS) part}
\begin{itemize}
\item $A$:  wet area
\item $Q$:  discharge
\item $\Omega(t,X)$:  free surface cross section
\item $H(t,X)$: physical water height
\item $h(t,X)$: $Z$-coordinate of the water level, $\sigma(X,h(t,X)) = T(A)$ : width of the free surface
\item $\nn_{\textbf{fs}}$: outward $\BB$-normal vector to the free surface
\item $\rho_0$:  density of the water at atmospheric pressure $p_0$
\end{itemize}
\subsubsection*{Notations concerning the pressurized part}
\begin{itemize}
\item $\Omega(X)$:  pressurized cross section
\item $\rho(t,X)$:  density of the water
\item $\beta$:  water compressibility coefficient
\item $c=\frac{1}{\sqrt{\beta\,\rho_0}}$:  sonic speed
\item $A = \frac{\rho}{\rho_0} S$:  FS equivalent wet area
\item $Q$:  FS equivalent discharge
\end{itemize}
\subsubsection*{Notations concerning the \textbf{PFS} model}
\begin{itemize}
\item $\SE$: the physical wet area: $\SE=A$ if the state is free surface, $S$ otherwise
\item $\mathcal{H}$: the $Z$ coordinate of the water level: $\mathcal{H}=h$ if the state is free surface, $R$ otherwise
\end{itemize}
\subsubsection*{Other notations}
\begin{itemize}
\item Bold characters are used for vectors except for $\SE$
\end{itemize}
}

\section{Formal derivation of the \textbf{FS}-model for free surface flows }\label{SectionFormalDerivationFSModel}
The classical shallow water equations are  used to describe physical situations like rivers, coastal domains, oceans and sedimentation problems.
These equations are obtained from  the  incompressible  Euler system (see e.g. \cite{AL08,LOE96}) or from the incompressible Navier-Stokes system (see for instance \cite{BCNV08,BN07,GP01,M07}) by several techniques (e.g. by direct integration or asymptotic analysis).  We adapt here the derivation in \cite{BFML07,BMCPV03} to get a  new unidirectional shallow water model. We start from  the $3$D incompressible Euler equations where we neglect the  acceleration following the $y$-axis supposing the existence of a privileged main flow axis. We write then the Euler equations in the local Serret-Frenet reference frame in order to take into account the local effects produced by the changes of section and the slope variation. Then we derive a shallow water model by a formal asymptotic analysis (done in Subsection \ref{SectionFSModelAsymptotic}).

\subsection{Incompressible Euler equations and framework}\label{SubsectionFramework}
Let us consider the cartesian refe\-rence frame $(O,\ii,\ji,\ki)$. In the corresponding co\-ordinate system  $(x,y,z)$, the  $3$D
incompressible Euler system writes:
\begin{equation}\label{E3DICartesian}
\left\{
\begin{array}{rcl} \dive(\rho_0\,\U) &=&0\\
\partial_t (\rho_0\,\U) + \rho_0\,\U\cdot\nabla (\rho_0\,\U) + \nabla P &=& \FF
\end{array}\right.
\end{equation}
where $\U(t,x,y,z)$ denotes the velocity with components $(u,v,w)$, $P=p(t,x,y,z)I_3$ is the isotropic pressure tensor, $\rho_0$ the density
of the fluid at atmospheric pressure $p_0$ and  $\FF$ is the exterior strength of gravity.

\noindent We close classically System~(\ref{E3DICartesian}) using a kinematic law for the evolution of the free surface: \emph{any free
surface particle is advected by the fluid velocity} $\U$ and on the wet boundary, we assume the no-leak condition $\U.\nn_{\textbf{wb}} = 0$
where $\nn_{\textbf{wb}}$ is the outward unit normal vector  to the wet part of the pipe (see \fig\ref{Supplementaire}).  We set the
pressure $P$ to $0$ at the free surface.

\begin{figure}[H]
 \begin{center}
 \includegraphics[scale = 0.6]{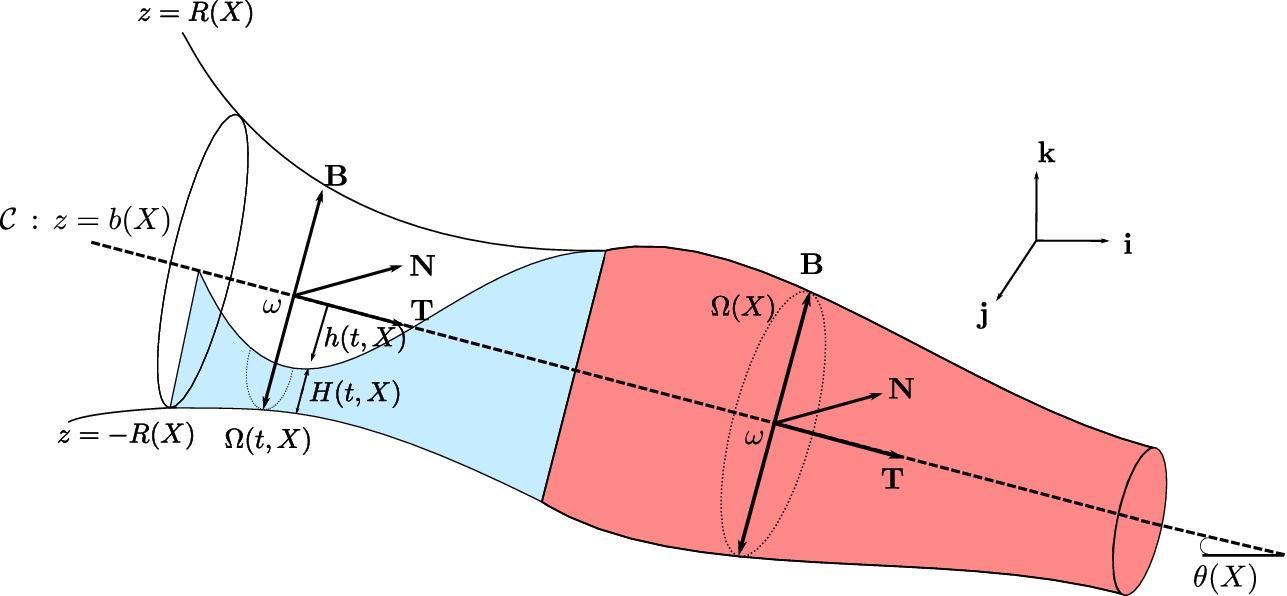}
  \caption{Geometric characteristics of the domain}
  \small{Mixed flow: free surface and pressurized}
  \label{OxOz}
 \end{center}
\end{figure}

\noindent We define the domain $\Omega_F(t)$  of the flow at time $t$ as the union of sections $\Omega(t,x)$
(assumed to be simply connected  compact sets) orthogonal  to some plane curve $\mathcal{C}$ lying in $(O,\ii,\ki)$ to follow the privileged
main flow axis. We  choose the parametrization
$(x,0,b(x))$  in the cartesian refe\-rence frame $(O,\ii,\ji,\ki)$ where $\ki$ follows the
vertical direction; $b(x)$ is then the elevation of the point $\omega (x,0,b(x))$ over the plane
$(O,\ii,\ji)$ (see \fig\ref{OxOz}).

\begin{figure}[H]
 \begin{center}
 \includegraphics[scale = 0.35]{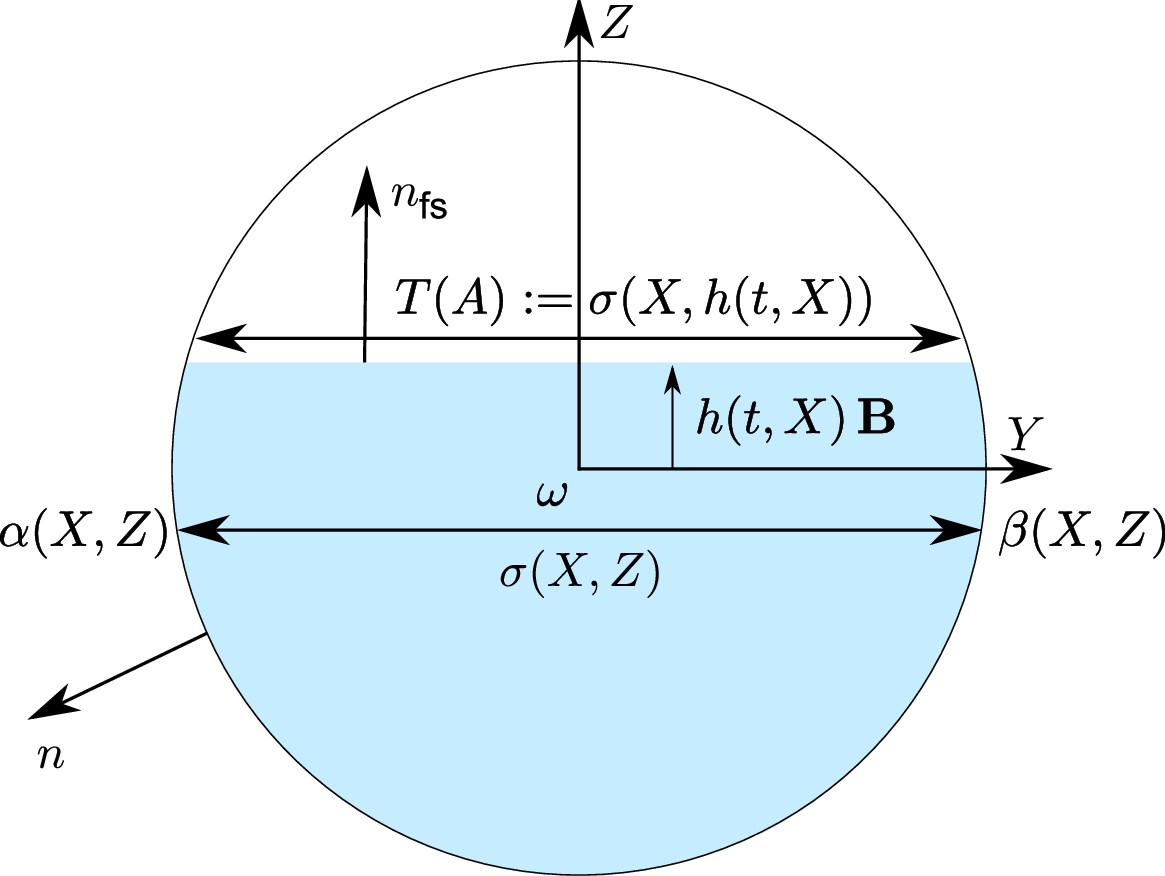}
  \caption{Cross-section $\Omega(t,X)$ of the domain at point $\omega$ in the free surface case}
  \label{OyOz}
 \end{center}
\end{figure}

\noindent We define a local reference frame as follows:  we introduce
the curvilinear variable defined by: $$\dsp X = \int_{x_0}^x \sqrt{1+(b'(\xi))^2} d\xi$$
where  $x_0$ is  an arbitrary abscissa. We set $Y=y$ and we denote by $Z$ the $\BB$-coordinate of any fluid
particle $M$ in the Serret-Frenet reference frame  $(\TT,\NN,\BB)$ at point $\omega (x,0,b(x))$ with  $\TT$ the tangent vector $\NN$, the
normal  and $\BB$ the binormal vector (see \fig\ref{OxOz} and \fig \ref{OyOz} for the notations). $\BB$ is normal to
$\mathcal{C}$ in the vertical plane $(O,\ii,\ki)$.

\noindent Then, at each point $\omega$, $\Omega(t,X)$ is defined by the set: $$\left\{(y,Z)\in\R^2; Z \in [-R(X),-R(X)+H(t,X)],\,
y \in [\alpha(X,Z),\beta(X,Z)]\right\}$$
\noindent where $R(X)$ denotes the radius, $H(t,X)$ the physical water height  at section $\Omega(t,X)$. We denote  $\alpha(X,Z)$ (respectively $\beta(X,Z))$ 
Y-coordinate of the left (respectively  right) boundary point of the domain at altitude $Z$, $-R(X)<Z<R(X)$ (see \fig\ref{OyOz}). We denote also
$-R(X)+H(t,X)$ by $h(t,X)$ which is the $Z$-coordinate of the water level.

\begin{figure}[H]
 \begin{center}
 \includegraphics[scale = 0.35]{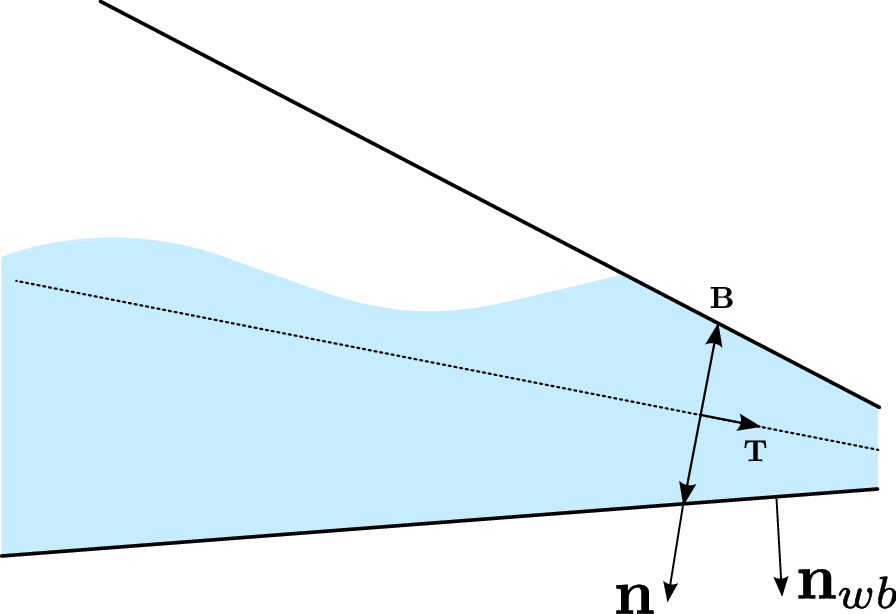}
  \caption{Outward unit normal $\nn_{\textbf{wb}}\neq \nn$ (except for uniform pipes)}
  \label{Supplementaire}
 \end{center}
\end{figure}

\noindent In the sequel, we will use a curvilinear map which will be an admissible transformation under the geometrical hypothesis
on the domain:
\begin{description}
\item[$(H)$] Let $\mathcal R(x)$ be the algebraic curvature radius of the plane curve $x\mapsto (x,0,b(x))$. We assume that:
$$\forall x \in \Omega_F,\,\,|\mathcal{R}(x)|> R(x).$$
\end{description}
\subsection{Incompressible Euler model in the curvilinear coordinates}\label{SubsectionEulerCurvilinear}
Following the work in \cite{BFML07,BMCPV03}, we write  System~(\ref{E3DICartesian}) in the Serret-Frenet reference frame  $(\TT,\NN,\BB)$ at point  $\omega(x,0,b(x))$ by the transformation $\mathcal{T} : (x,y,z)\rightarrow (X,Y,Z)$ using the divergence chain rule lemma that we recall here:
\begin{lem}\label{lemma}

Let $(X,Y,Z) \mapsto \mathcal{T}(X,Y,Z)=(x,y,z)$ be a $C^1$  diffeomorphism and \\$\mathcal{A}^{-1} = \nabla_{(X,Y,Z)} \mathcal T$ the jacobian
matrix of the transformation with determinant $J$.

\noindent Then, for any vector field $\Phi$, one has:
$$J \dive_{(x,y,z)} \Phi = \dive_{(X,Y,Z)}(J\mathcal{A}\Phi)\,,$$
and, for any scalar function $f$, one has:
$$\nabla_{(x,y,z)}f = \mathcal{A}^t \nabla_{(X,Y,Z)}f,$$ where $\mathcal{A}^t$ stands for the transpose of the matrix $\mathcal{A}$.
\end{lem}

\noindent Let $(U,V,W)^t$ be the components of the velocity vector  in the $(X,Y,Z)$ coordinates  defined as $(U,V,W)^t = \Uptheta (u,v,w)^t$ where  $\Uptheta$ is the matrix $$\Uptheta=\left(
\begin{array}{ccc}
 \cos\theta & 0 & \sin\theta \\
 0 & 1 & 0 \\
 -\sin\theta& 0 & \cos\theta \\
\end{array} \right)\,,$$  where we  denote by $\theta(x)$ the angle $(\ii,\TT)$ in the $(\ii,\ki)$ plane.

\noindent Using Lemma \ref{lemma}, the incompressible  Euler system  in the variables $(X,Y,Z)$ reads:
\begin{equation}\label{E3DICurvilinear}\left\{
\begin{array}{lll}
\partial_X (\rho_0\,U) + \partial_Y(J\rho_0\,V)+  \partial_Z(J\rho_0\,W)&=&0\\
\partial_t(J\rho_0\,U) + \partial_X(\rho_0\,U^2) + \partial_Y(J\rho_0\,UV)+  \partial_Z(J\rho_0\,UW)+ \partial_X p &=& G_1\\
\partial_t(J\rho_0\,V) + \partial_X(\rho_0\,UV) + \partial_Y(J\rho_0\,V^2)+  \partial_Z(J\rho_0\,VW)+ \partial_Y(Jp) &=& 0\\
\partial_t(J\rho_0\,W) + \partial_X(\rho_0\,UW) + \partial_Y(J\rho_0\,VW)+  \partial_Z(J\rho_0\,W^2)+ J\partial_Z(p) &=& G_2
\end{array}\right.\end{equation} where $J(X,Y,Z) = 1-Z\dsp \theta'(X)$ is the determinant of the transformation and $$G_1 =
\rho_0\,UW\theta'(X) -Jg\rho_0\,\sin\theta,\,\,  G_2=-\rho_0\,U^2\theta'(X)-Jg\rho_0\,\cos\theta.$$ The interested reader can find
the details of the calculus in \cite{BFML07}. We have denoted by $f'$ the  derivative with respect to the space variable $X$ of any function
$f(X)$.

\noindent On the wet boundary, the no-leak condition reads:
\begin{eqnarray}
({U},{V},{W})^t.\nn_\textbf{wb}=0\,. \label{NoLeakCond}
\end{eqnarray}

\begin{rque}\label{RemarkAdmissibilityCondition}

\noindent Notice that $\kappa(X)=\dsp\theta'(X)$ is the algebraic curvature of the axis at point $\omega(X,0,b(X))$ and the
function $J(X,Y,Z) = 1-Z \kappa(X)$ depends only on the variables $X,Z$. Moreover, under the hypothesis $(H)$, we have $J>0$ in $\Omega_F$.
Consequently, $\mathcal T$ defines a diffeomorphism and thus the performed transformation is admissible.
\end{rque}

\subsection{Formal derivation of the \textbf{FS}-model for free surface  flows}\label{SectionFSModelAsymptotic}
In this section, we perform a formal asymptotic analysis on System~(\ref{E3DICurvilinear}). According to the work in \cite{BFML07,GP01,M07}, the shallow water equations  can be obtained from the incompressible Navier-Stokes equations with particular boundary conditions. Here, we perform this analysis directly on the incompressible Euler system in order to get $J = 1 + O(\epsilon)$ for some small parameter $\epsilon$.

\noindent Let us introduce the usual small parameter $\epsilon = H/L$ where $H$ (the height) and $L$ (the length) are two characteristics
dimensions along the $\BB$ and $\TT$ axis respectively. Moreover, we  assume that the characteristic dimension along the $\ji$ axis is the same as for the $\ki$ axis. 
We introduce the other characteristics dimensions $T,P,\overline U,\overline V,\overline W$ for time, pressure and velocity
respectively and the dimensionless quantities as follows: $$\widetilde U = U/ \overline U,\, \widetilde V = \epsilon V/ \overline U,\,
\widetilde W = \epsilon {W}/{\overline U},\, $$ $$\widetilde X =
{X}/{L},\, \widetilde Y = {Y}/{H},\, \widetilde Z = {Z}/{H},\,
\widetilde p = {p}/{P},\widetilde\theta = {\theta},\widetilde\rho = {\rho_0}.$$ In the sequel, we set $P={\overline
U}^2$ and $L=T\overline U$ (i.e. we only consider  laminar flows).

\noindent Under these hypotheses, we have $\widetilde J (\widetilde X,\widetilde Y,\widetilde Z)
= 1- \epsilon \widetilde Z\dsp \widetilde\theta'(\widetilde X) $. Thus, the rescaled System~(\ref{E3DICurvilinear}) reads:
\begin{equation}\label{E3DICurvilinearRescaled}\left\{
\begin{array}{rcl}
\partial_{\widetilde X} \widetilde U + \partial_{\widetilde Y}(\widetilde J \widetilde V)+  \partial_{\widetilde Z}
(\widetilde J \widetilde W)&=&0\\
\partial_{\widetilde t}(\widetilde J \widetilde U) + \partial_{\widetilde X}({\widetilde U}^2) +
\partial_{\widetilde Y}(\widetilde J \widetilde U \widetilde V)+
\partial_{\widetilde Z}(\widetilde J \widetilde U \widetilde W)+ \partial_{\widetilde X} \widetilde p &=& G_{1}\\
\epsilon^2\left(\partial_{\widetilde t}(\widetilde J \widetilde V) +
\partial_{\widetilde X}(\widetilde U \widetilde V) + \partial_{\widetilde Y}
(\widetilde J {\widetilde V}^2)+  \partial_{\widetilde Z}(\widetilde J \widetilde V \widetilde W)\right)+
\partial_{\widetilde Y}(\widetilde J \widetilde p) &=& 0\\
\epsilon^2\left(\partial_{\widetilde t}(\widetilde J \widetilde W) +
\partial_{\widetilde X}(\widetilde U \widetilde W) + \partial_{\widetilde
Y}(\widetilde J \widetilde V \widetilde W) +  \partial_{\widetilde Z}(\widetilde J
{\widetilde W}^2)\right)\\  + \widetilde J\partial_{\widetilde Z}(\widetilde p) = G_{2}
\end{array}\right.
\end{equation}
\noindent where $$G_{1} =  \epsilon \widetilde U \widetilde W \widetilde \kappa(\widetilde
  X)-\dsp\frac{\sin\widetilde\theta}{{F_{r,L}}^2}-
  \frac{\widetilde Z \dsp }{{F_{r,H}}^2} (\cos\widetilde\theta)', $$
$$G_{2} =
-\epsilon{\widetilde U}^2\widetilde\rho(\widetilde X)-
\dsp\frac{\cos\widetilde\theta}{{F_{r,H}}^2}
+
\dsp\epsilon\kappa(X)\frac{\widetilde Z\widetilde J\cos\widetilde\theta}{{F_{r,H}}^2},$$
 $F_{r,M} = \dsp\frac{\overline{U}}{\sqrt{gM}}$ is the Froude number along the $\TT$ axis and the $\BB$ or $\NN$ axis where $M$ is
any generic variable equal to $L$  or $H$.

\noindent Formally, when $\epsilon $ vanishes, System (\ref{E3DICurvilinearRescaled}) reduces to:
\begin{eqnarray}\label{E3DICurvilinearRescaledEpsilonVanishes}
\partial_{\widetilde X} \widetilde U + \partial_{\widetilde Y}(\widetilde V)+  \partial_{\widetilde Z}(\widetilde W)&=&0
\label{E3DICurvilinearRescaledEpsilonVanishes1} \\
\partial_{\widetilde t}(\widetilde U) + \partial_{\widetilde X}(\widetilde U^2) + \partial_{\widetilde Y}(\widetilde U
\widetilde V)+ \partial_{\widetilde Z}(\widetilde U
\widetilde W)+
\partial_{\widetilde X} \widetilde p &=&
  -\dsp\frac{\sin\widetilde\theta}{{F_{r,L}}^2} \nonumber\\
 & & - \dsp  \frac{\widetilde Z \dsp}{{F_{r,H}}^2} (\cos\widetilde\theta)'
\label{E3DICurvilinearRescaledEpsilonVanishes2}\\
\partial_{\widetilde Z}(\widetilde p) &=&
-\dsp\frac{\cos\widetilde\theta}{{F_{r,H}}^2} \label{E3DICurvilinearRescaledEpsilonVanishes3}
\end{eqnarray}

\noindent Let us introduce the  conservative variables $A(t,X)$ and $Q(t,X)$  representing respectively the wet area and the discharge  defined as: $$A(t,X) = \int_{\Omega(t,X)} dYdZ,\quad Q(t,X)=A(t,X)\overline{U}$$
where  $\overline{U}$ is the mean value of the velocity :$$\overline{U}(t,X)=\frac{1}{A(t,X)}\int_{\Omega(t,X)}U(t,X)\,\,dY dZ\,.$$
 We integrate the preceding system
(\ref{E3DICurvilinearRescaledEpsilonVanishes1}-\ref{E3DICurvilinearRescaledEpsilonVanishes2}-\ref{E3DICurvilinearRescaledEpsilonVanishes3})
along the cross-section with the approxi\-mation $\overline{U^2}\approx \overline{U}\,\overline{U}$ and $\overline{U\,V}\approx
\overline{U}\,\overline{V}$. Then, returning to the physical variables, the free surface model, that we call \textbf{FS}-model, reads:
\begin{equation}\label{FSModel}
\left\{\begin{array}{lll}
\partial_{t}A +\partial_X Q &=& 0\\
\dsp\partial_t Q +\dsp\partial_X\left(\frac{Q^2}{ A}+\dsp g I_1(X,A) \cos\theta  \right)&=&\dsp g I_2(X,A) \cos\theta
-gA\sin\theta \\ & &-g A \overline{Z}(X,A) \dsp{(\cos\theta)'}\end{array}\right.
\end{equation}
where  $I_1(X,A)$ and $I_2(X,A)$ are respectively the classical term of hydrostatic pressure and the pressure source term  defined by:  $$I_1(X,A) =
\int_{-R}^{h}(h-Z) \sigma \,dZ \textrm{ and } I_2(X,A) =
\int_{-R}^{h} (h-Z) \partial_X\sigma \,dZ$$ which are obtained from the integration of the pressure term
$\partial_{\widetilde{X}}\widetilde{p}$ in Equation (\ref{E3DICurvilinearRescaledEpsilonVanishes1}) with  $\widetilde{p} = \rho(h(t,X)-Z)\cos\theta$ (obtained from equation
(\ref{E3DICurvilinearRescaledEpsilonVanishes3})).

\noindent In these formulas $\sigma(X,Z)$ is the width of the cross-section at position $X$ and at height $Z$. The additional  term
$\overline{Z}(X,A)$ is defined by $(h(A)-I_1(X,A)/A)$. It is the $Z$-coordinate of the center of mass:
$$
\begin{array}{lll}
\overline{Z} &=& \dsp\int_{\Omega(t,X)} Z \,dY\,dZ\\
             &=& \dsp\int_{-R(X)}^{h(t,X)} \int_{\alpha(X,Z)}^{\beta(X,Z)} Z \,dY\,dZ\\
             &=& \dsp\int_{-R(X)}^{h(t,X)} Z\, \sigma(X,Z)\,dZ
\end{array}\,.
$$
In System (\ref{FSModel}), we may add a friction term $-\rho_0 g S_f\, \TT$ to take into account the 
dissipation of energy. We have chosen this term $S_f$ as the one  given by the Manning-Strickler law (see e.g. \cite{SWB98}):
 $$S_f(A,U)=K(A)U|U|\,.$$ 
 The term $K(A)$ is defined by: $\dsp K(A) = \frac{1}{K_s^2 R_h(A)^{4/3}}$, $K_s>0$ is the Strickler coefficient of  roughness depending on the material, 
  $R_h(A)= A/P_m$ is the hydraulic radius and $P_m$ is the perimeter of the wet surface area 
  (length of the part of the channel's section in contact with the water).

\section{Formal derivation of the \textbf{P}-model for pressurized flows}\label{SectionPFModel}
In this section, we present a new set of unidirectional shallow water like equations to describe pressurized flows in closed non uniform water pipes. This model is constructed to be coupled in natural way with the \textbf{FS}-model~(\ref{FSModel}). Starting from the $3$D compressible Euler equations in cartesian coordinates,
\begin{equation}\label{Euler3D_mass_conservation}
\partial_t \rho + \dive{(\rho \U)} = 0,
\end{equation}
\begin{equation}\label{Euler3D_momentum_conservation}
\partial_t (\rho\U) + \dive{(\rho \U\otimes \U)} + \nabla p = \FF,
\end{equation}
where $\U(t,x,y,z)$ and $\rho(t,x,y,z))$ denotes the velocity with components $(u,v,w)$ and  the density respectively. $p(t,x,y,z)$ is the scalar pressure and $\FF$ the exterior strength of gravity.

\noindent We define the pressurized domain of the flow as the continuous extension of $\Omega_F$ (see Subsection~\ref{SubsectionFramework}) defined by some plane curve $\mathcal{C}$ with parametrization $(x,0,b(x))$ in the cartesian reference frame $(O,\ii,\ji,\ki)$; we recall that $b(x)$ is then the elevation of the point $\omega$ over the plane $(O,\ii,\ji)$ (see \fig\ref{OxOz}). The curve may be, for instance, the axis spanned by the center of mass of each orthogonal section $\Omega(x)$ to the main mean flow axis, particularly in the case of a piecewise cone-shaped pipe. Notice that we consider only the case of infinitely rigid pipes, thus the sections $\Omega=\Omega(x)$ are only $x$-dependent.

\noindent We then write Equations (\ref{Euler3D_mass_conservation}-\ref{Euler3D_momentum_conservation}) in the $(X,Y,Z)$ coordinates introduced in Subsection~\ref{SubsectionFramework}. As we want a unidirectional model, we suppose that the mean flow follows the $X$-axis. To this end, we neglect the second and third equation for the conservation of the momentum.\\
\noindent By a straightforward computation, the mass and the first momentum conservation equation in the $(X,Y,Z)$ coordinates reads:

\begin{equation}\left\{
\begin{array}{rcl}
\partial_t (J\rho) + \partial_X(\rho U) + \partial_Y(\rho JV)+ \partial_Z(\rho JW) &=& 0\\
& & \label{Euler_system_curvilinear}\\
\partial_t(J\rho U) + \partial_X(\rho U^2) + \partial_Y(\rho JUV^2)+ \partial_Z(\rho JUW)+ \partial_X p \\ = -\rho Jg\sin\theta  +\rho U W  \dsp (\cos\theta)'\end{array}\right.\end{equation}

\noindent Applying the same asymptotic analysis developed in Subsection \ref{SectionFSModelAsymptotic}, Equations (\ref{Euler3D_mass_conservation}-\ref{Euler3D_momentum_conservation}) read:

\begin{equation}\label{E3DCRescaledCurvilinear}\left\{
\begin{array}{rcl}
\partial_{t}(\rho) + \partial_{X}({\rho U}) +
\partial_{Y}(\rho V)+
\partial_{Z}(\rho W)&=&0\\
& & \\
\partial_{ t}(\rho U  ) + \partial_{  X}({ \rho  U}^2) +
\partial_{  Y}( \rho  U   V)+
\partial_{  Z}( \rho  U   W)+ \partial_{  X}   p &=&
  -\rho\dsp{g\sin \theta} \\ & & -
   g Z  \dsp (\cos\theta)'
\end{array}\right.\end{equation}

\noindent We choose the linearized pressure law:
\begin{equation}\label{LinearizedPressureLaw}
p = p_a + \dsp{\frac{\rho-\rho_0}{\beta \rho_0}}
\end{equation}
(see e.g. \cite{SWB98,WS78}) in which $\rho_0$ represents the
density of the fluid at atmospheric pressure $p_0$, $p_a$ is some function set to zero and $\beta$ the water compressibility coefficient  (equal to $5.0\,10^{-10}\,m^2.N^{-1}$
in practice). The sonic speed is then given by $c~=~1/\sqrt{\beta \rho_0}$ and thus $c \approx 1400\,m.s^{-1}$.\\
For $m \in \partial \Omega$, $\dsp \nn = \frac{\mm}{|\mm|}$ is the outward unit vector at the point $m$ in the $\Omega$-plane and $\mm $
stands for the vector \textbf{${\omega m}$} (as displayed on \fig\ref{OyOz}).

\noindent Following the section-averaging method performed in Subsection \ref{SectionFSModelAsymptotic}, we integrate
System~(\ref{E3DCRescaledCurvilinear})  over the cross-section $\Omega$. Noting  the averaged values over $\Omega$ by the overlined letters
(except $\overline{Z}$), and using the approximations $\overline{\rho U}\approx \overline{\rho}\overline{U},\,\overline{\rho U^2}\approx
\overline{\rho}\overline{U}^2$ the shallow water like equations read:
\begin{eqnarray}
\partial_{t}(\overline\rho S) + \partial_{X}({\overline\rho S \overline U}) & = & \dsp \int_{\partial \Omega} \rho \left( U\partial_X \mm - \VV\right).\nn\, ds \label{STV_mass}\\
\partial_{t}(\overline\rho S \overline U) + \partial_{X}(\dsp\overline\rho S \overline{U}+c^2\overline\rho S) &=&
-\dsp{g\overline\rho S \sin \theta}+c^2\overline\rho S'\nonumber\\	
 & -&   g\overline\rho S\overline Z  \dsp (\cos\theta)'\label{STV_momentum}\\
 &  +&   \dsp \int_{\partial \Omega} \rho U \left(U\partial_X \mm - \VV\right).\nn\, ds\nonumber
\end{eqnarray}
where $\VV = (V,W)^t$ is the velocity  in the $(\NN, \BB)$-plane. We denote by $S$ the area of the cross-section  of the pipe at position
$X$.

\noindent The integral terms appearing in (\ref{STV_mass}) and (\ref{STV_momentum}) vanish, as the pipe is infinitely rigid, i.e. $\Omega = \Omega(X)$ (see \cite{BG08} for the dilatable case). It follows the  non-penetration condition (see \fig\ref{Restriction}): $$\vectrois{U}{V}{W}.\nn_{\textbf{wb}} = 0 \,.$$

\noindent Omitting the overlined letters (except $\overline{Z}$), we introduce the conservative variables
\begin{eqnarray}
A = \dsp\frac{\rho }{\rho_0} S& \textrm{the FS \emph{equivalent wet area}} \label{EquivalentWetArea}\\
Q =\dsp A U & \textrm{ the FS \emph{equivalent discharge}} \label{EquivalentWetDischarge}\,.
\end{eqnarray}
and dividing Equations (\ref{STV_mass})-(\ref{STV_momentum}) by $\rho_0$ we get:
\begin{equation}\label{PModel}\left\{
\begin{array}{rcl}
\partial_{t}(A) + \partial_{X}(Q) &=&0\\
&  & \\
\partial_{t}(Q) + \partial_{X}\left(\dsp\frac{Q^2}{A}+c^2 A\right) &=&
  -\dsp{g A \sin\theta}-
gA \overline{Z}(X,S)   \dsp (\cos\theta)'  \\ &  & +\dsp c^2 A \frac{S'}{S}
\end{array}\right.\end{equation}
As introduced previously for the \textbf{FS}-model in Section (\ref{SectionFSModelAsymptotic}), we may introduce the friction term  $-\rho g S_f\, \TT$ given by the Manning-Strickler law (see e.g. \cite{SWB98}): $$S_f(S,U)=K(S)U|U|$$ where $K(S)$ is defined by: $\dsp K(S) = \frac{1}{K_s^2 R_h(S)^{4/3}}$, $K_s>0$ is the Strickler coefficient
of roughness depending on the material and $R_h(S)= S/P_m$ is the hydraulic radius where $P_m$ is the perimeter of the wet surface area (length of the part of the channel's section in contact with the water, equal to $2\,\pi\,R$ in the case of circular pipe).

\paragraph*{}
This choice of variables is motivated by the fact that this system is formally close to the
\textbf{FS}-model (\ref{FSModel}) where the terms $gI_1(X,A) \cos\theta$, $gI_2(X,A) \cos\theta$,
$\overline{Z}(X,A)$ are respectively the counterparts of $c^2 A$, $\dsp c^2 A
\frac{S'}{S}$, $\overline{Z}(X,S)$ in System~(\ref{PModel}). Let us remark that the term $\overline{Z}$ is continuous through the change of
state (pressurized to  free surface or free surface to pressurized state) when the same curve plane is chosen (in practice, the main axis
of the pipe). Then, we are motivated to connect ``continuously'' System (\ref{FSModel}) and (\ref{PModel}) through transition points (through the change of state) by defining a continuous pressure law. It leads to a ``natural'' coupling between the pressurized and free
surface model as we will see in Section~\ref{SectionCoupling}.

\section{The \textbf{PFS}-model}\label{SectionCoupling}
The formulations of the \textbf{FS}-model~(\ref{FSModel}) and \textbf{P}-model~(\ref{PModel}) are very close to each other. The main difference comes from the
pressure law. In order to build a coupling between the two models, we have to define a pressure that ensures its continuity through
transition points in the same spirit of \cite{BG07}. As pointed out in the previous section, we will use the common couple of unknowns $(A,Q)$
and the same plane curve  $\mathcal{C}$ (see Remark \ref{RemarkRestrictionPlaneCurve}) to get a continuous model for mixed flows.
\begin{rque}\label{RemarkRestrictionPlaneCurve}
\noindent The plane curve with parametrization $(x,0,b(x))$ is chosen as the main pipe axis in the axisymmetric case. Actually this choice is the more convenient for pressurized flows while the
bottom line is adapted to free surface flows. Thus we must assume small variations of the section
($\dsp S'$ small) or equivalently small angle $\varphi$ as displayed on \fig\ref{Restriction}.

\begin{figure}[H]
 \begin{center}
 \includegraphics[scale = 0.6]{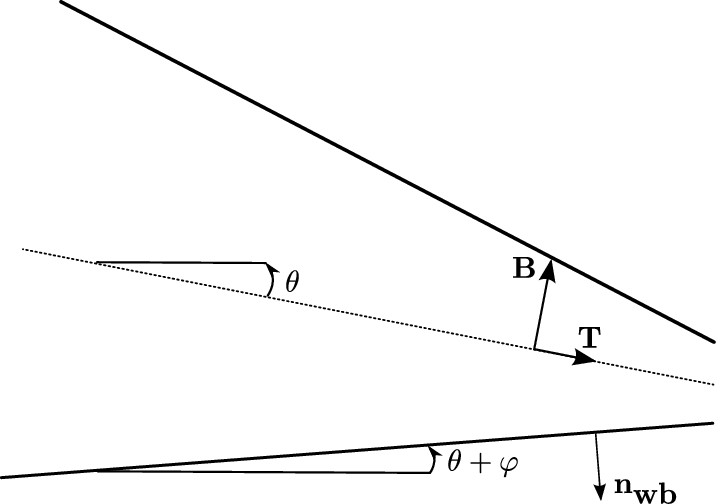}
  \caption{Some restriction concerning the geometric domain.}
  \label{Restriction}
 \end{center}
\end{figure}
\end{rque}

\begin{figure}[H] \label{Depression}
\centering
\includegraphics[angle=0,scale=0.45]{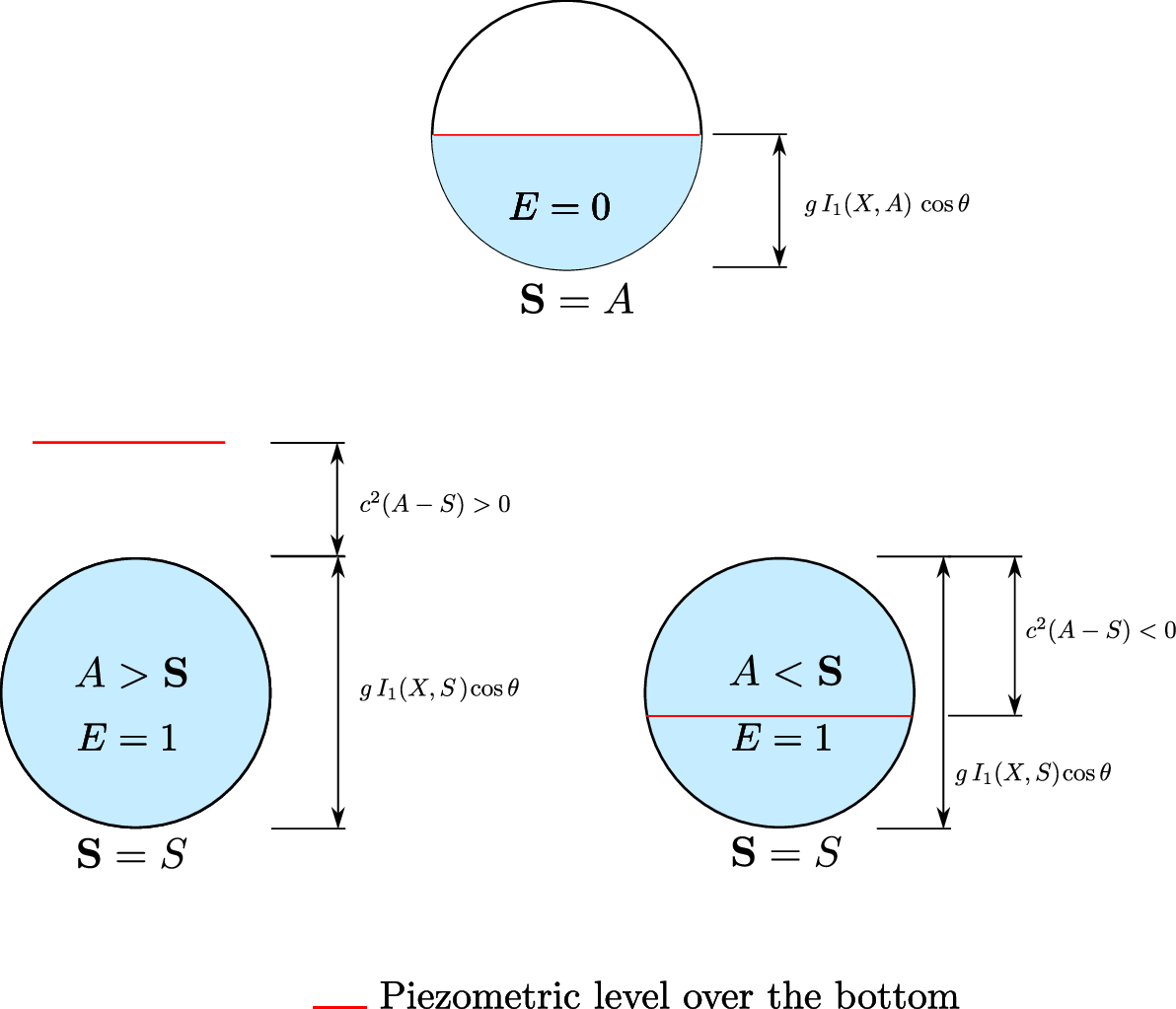}
\caption{Free surface state $p(X,A,0) = g\,I_1(X,A)\cos\theta$ (top), pressurized state with overpressure $p(x,A,1)>0$ (bottom left),
pressurized state with  depression
$p(x,A,1)<0$ (bottom right).}
\end{figure}

\noindent We introduce a state indicator $E$ (see \fig\ref{Depression}) such that:
\begin{equation}\label{DefinitionE}
E =\left\{
\begin{array}{ll}
1 & \textrm{ if the state is pressurized: } (\rho\neq\rho_0) \\
0 & \textrm{ if the state is free surface: } (\rho=\rho_0)\\
\end{array}
\right.\,.
\end{equation}

Next, we define the \emph{physical wet area} $\SE$ by:
\begin{equation}\label{DefinitionPhysicalWetArea}
\SE = \SE(A,E) = \left\{
\begin{array}{lll}
S       & \textrm{ if } & E = 1 \\
A       & \textrm{ if } & E = 0
\end{array}
\right.
\end{equation}
and a modified pressure law (see \fig\ref{Depression}) which ensures its  continuity  through the change of state by:
\begin{equation}\label{PFSModelPressureInFlux}
 p(X,A,E) =  c^2(A-\SE) + gI_1(X,\SE)\cos\theta .
\end{equation}

\begin{rque}
\begin{itemize}
\item[]
\item Indeed, when a change of state occurs we have:
$$\lim_{\stackrel{A\to S}{A<S}} p(X,A,E) = \lim_{\stackrel{A\to S}{A>S}} p(X,A,E) = g I_{1}(X,S) \cos(\theta) $$ 
which ensures the continuity
of the pressure.
\item The flux gradient $F$ is discontinuous through
the change of state since $$\frac{\partial F}{\partial A}(A,Q,0) = g\frac{\partial}{\partial A}I_1(X,A)\cos\theta \neq c^2=\frac{\partial F}{\partial A}(A,Q,1) .$$
\end{itemize}
\end{rque}

\noindent Finally, from the \textbf{P}-model (\ref{PModel}), the \textbf{FS}-model (\ref{FSModel}), 
the definition of $E$ (\ref{DefinitionE}), the definition of $\SE$
 (\ref{DefinitionPhysicalWetArea}) and the pressure law (\ref{PFSModelPressureInFlux}),  the \textbf{PFS}-model for unsteady mixed flows can be
simply expressed into a single formulation as:
\begin{equation}\label{PFS}
\left\{
\begin{array}{rcl}
\partial_{t}(A) + \partial_{X}({Q}) &=&0\\
\partial_{t}(Q) + \partial_{X}\left(\dsp\frac{Q^2}{A}+p(X,A,E)\right) &=&
-\dsp g A\dsp b' + Pr(X,A,E)\\ & & -G(X,A,E)\\ & &- K(X,A,E)\dsp\frac{Q|Q|}{A}
\end{array}\right.
\end{equation}
\noindent where $K$, $Pr$, and $G$ denotes respectively the friction, the pressure source  and the geometry source term defined as
follows:
$$\begin{array}{lll}
Pr(X,A,E) &=& \dsp c^2\left(\frac{A}{\SE}-1\right)\, S'  + g I_2(X,\SE) \cos\theta \\
          & & \textrm{ with } I_2(X,\SE) = \dsp\int_{-R(X)}^{\mathcal{H}(\SE)}
(\mathcal{H}(\SE)-Z)\,\partial_X \sigma(X,Z)\,dZ,\\
 G(X,A,E) &=& \dsp gA\,  \overline{Z}(X,\SE) \dsp (\cos\theta)', \\
K(X,A,E) &=& \dsp \frac{1}{K_s^{2} R_h(\SE)^{4/3}}\,
  \end{array}
$$
and $b'$ stands for $\sin\theta(X)$. $\mathcal{H}$ represents  the $Z$-coordinate of the water level:
\begin{equation}\label{DefinitionPhysicalWaterHeight}
\mathcal{H} = \mathcal{H}(\SE)=\left\{
\begin{array}{lll}
h(A)       & \textrm{ if } & E = 0 \\
R(X)       & \textrm{ if } & E = 1
\end{array}
\right.\,.
\end{equation}

\begin{rque}[Both models are recovered]\label{RqueRecoveringBothModels}

\noindent Setting $\SE(A,E) = A$ in System (\ref{PFS}), we obtain obviously the free surface model  (\ref{FSModel}).
For all pressurized states, when $\SE(A,E) = S$, the pressure law (\ref{PFSModelPressureInFlux}) reads, for instance, in the case of circular pipe: $$c^2(A-S)
+
gI_1(X,S)\cos\theta = c^2(A-S) + g\,\pi\,R^3\,\cos\theta$$ which is not exactly the pressure law of the \textbf{P}-model
(\ref{PModel}). Indeed, the derivation of the \textbf{P}-model is done with the linearized pressure law (\ref{LinearizedPressureLaw}) (see
Section \ref{SectionPFModel}) with $p_a=0$. Thus, the property of the continuity of models (\ref{PModel})-(\ref{FSModel}) through a
change of state is obtained if and only if $p_a$ is chosen as $gI_1(X,S)\cos\theta$ which is the hydrostatic pressure corresponding to a full section.
\end{rque}

\newpage
\noindent The \textbf{PFS}-model (\ref{PFS}) satisfies the following properties:
\begin{thm}\label{ThmPFSModel}
\begin{enumerate}
\item[]
\item The right eigenvalues of System (\ref{PFS}) are given by:
$$\lambda^- =U-c(A,E),\,\lambda^+ = U + c(A,E)
$$
with $c(A,E)= \left\{
\begin{array}{lll}
\dsp\sqrt{g\,\frac{A}{T(A)}\,\cos\theta} &\textrm{ if } & E=0\\
\dsp c  &\textrm{ if } & E=1
\end{array}
\right.\,,
$
where $T(A)$ is the width of the free surface (see \fig\ref{OyOz}).

\noindent Then, System (\ref{PFS}) is strictly hyperbolic on the set:
$$\left\{A(t,X)>0\right\}\,.$$
\item For smooth solutions, the mean velocity $U = Q/A$ satisfies
\begin{equation}\label{ThmPFSEquationForU}
\begin{array}{c}
\partial_t U + \partial_X \left(\dsp\frac{ U^2}{2} + c^2 \ln(A/\SE) + g\mathcal{H}(\SE)\cos\theta + gb \right) \\ =
-g K(X,A,E) U|U| \leqslant 0.
\end{array}
\end{equation}
\noindent The quantity $\dsp \frac{ U^2}{2} + c^2 \ln(A/\SE) + g\mathcal{H}(\SE)\cos\theta + gb$ is called the
total head.
\item The still water steady state reads:
\begin{equation}\label{ThmPFSSteadyState}
u = 0 \;\mbox{ and } \; c^2 \ln(A/\SE) + g\mathcal{H}(\SE)\cos\theta + gb = 0.
\end{equation}
\item It admits a mathematical entropy
\begin{equation}\label{ThmPFSMathematicalEntropy}
\mathcal{E}(A,Q,E) =\dsp \frac{Q^2}{2A} + c^2 A \ln(A/\SE)+ c^2 S + g A
\overline{Z}(X,\SE)\cos\theta + gAb
\end{equation}
\noindent which satisfies the entropy relation for smooth solutions
\begin{equation}\label{ThmPFSEntropy}
\partial_t \mathcal{E} +\partial_X \Big(\left(\mathcal{E}+p(X,A,E)\right)U\Big) = -gAK(X,A,E) U^2 |U| \leqslant 0\,.
\end{equation}
\end{enumerate}
\end{thm}
Notice that the total head and $\mathcal{E}$ are defined continuously through the transition points.

\begin{rque}
The  term $A\overline{Z}(X,A) (\cos\theta)'$ is also called ``corrective term'' since it allows to write the  Equations
(\ref{ThmPFSEquationForU}) and (\ref{ThmPFSEntropy}) with (\ref{ThmPFSMathematicalEntropy}).
\end{rque}

\noindent \textbf{Proof of Theorem~\ref{ThmPFSModel}:} the results~(\ref{ThmPFSEquationForU}) and (\ref{ThmPFSEntropy})  are obtained in a
classical way. Indeed, Equation (\ref{ThmPFSEquationForU}) is obtained by subtracting the result of the multiplication of the mass  equation
by $U$ to the momentum equation. Then multiplying the mass  equation by $\left(\dsp\frac{ U^2}{2} + c^2 \ln(A/\SE) +
g\mathcal{H}(\SE)\cos\theta + gb \right)$ and adding the result of the multiplication of Equation (\ref{ThmPFSEquationForU}) by $Q$,
we get:
$$\begin{array}{l}
\partial_t \left(\dsp \frac{Q^2}{2A} +c^2 A \ln(A/\SE)+ c^2 S +g A \overline{Z}(X,\SE)\cos\theta + gAb\right) \\
+\partial_X
\left(\left(\dsp \frac{Q^2}{2A} + c^2 A \ln(A/\SE)+ c^2 S +gA \overline{Z}(X,\SE)\cos\theta +
gAb+p(X,A,E)\right)U\right) \\
+c^2\left(\dsp\frac{A}{\SE}-1\right)\partial_t \SE = -gAK(X,A,E) U^2 |U| \leqslant 0 \,.\end{array}$$ We see that
the
term  $c^2\left(\dsp\frac{A}{\SE}-1\right)\partial_t \SE$ is identically $0$ since we have  $\SE=A$   when
the flow is free
surface whereas $\SE=S(X)$ when the flow is pressurized. Moreover, from the last inequality, when $\SE=A$, we have the
classical entropy
inequality (see \cite{BG07,BG08}) with $\mathcal{E}$: $$\mathcal{E}(A,Q,E) =\dsp \frac{Q^2}{2A} +gA \overline{Z}(X,A)\cos\theta +
gAb $$
while in the pressurized case, it is: $$\mathcal{E}(A,Q,E) =\dsp \frac{Q^2}{2A} + c^2 A \ln(A/S)+ c^2 S  + gAb.$$
Finally, the entropy for the \textbf{PFS}-model reads: $$\mathcal{E}(A,Q,E) =\dsp \frac{Q^2}{2A} + c^2 A \ln(A/\SE)+ c^2 S  +gA
\overline{Z}(X,\SE)\cos\theta + gAb.$$
Let us remark that the term $c^2 S$ makes $\mathcal{E}$ continuous through transition points and it permits also to write   the
entropy flux under the classical form $(\mathcal{E}+p)U$.
 \begin{flushright}
  $\blacksquare$
 \end{flushright}

\section{Finite volume discretisation}\label{PFSRoe}
In this section, we adapt the VFRoe scheme described in \cite{BGH00,GHS01,BG07}. The new terms appearing in the \textbf{PFS}-model related to the
curvature and the section variation are  upwinded in the same spirit of \cite{BG07}. The numerical scheme is adapted to disconti\-nui\-ties of
the flux gradient occurring in the treatment of the transitions between free surface and pressurized states.
\subsection{Discretisation of the space domain}
The spatial domain is a pipe of length $L$. The main axis of the pipe is divided in  cells $\dsp
m_i=[X_{i-1/2},X_{i+\frac{1}{2}}],\ 1\leq i\leq N$. $\Delta t^n$ denotes the timestep and we set
$t_{n+1}=t_n + \Delta t^n$.\\
The discrete unknowns are $U_i^n=\left(\begin{array}{c}A_i^n\\
Q_i^n\end{array}\right)$. For the sake of simplicity, the boundary conditions are not treated (the interested reader can find this treatment in details in \cite{BG07}).
%
\subsection{Explicit first order VFRoe scheme}\label{SubSectionExplicitFirstOrderRoeScheme}
We propose to extend the finite volume discretisation \cite{BG07} to the \textbf{PFS}-model using the upwinding of the new source terms: the curvature and section variation of the pipe. In what follows, we do not write the $E$ dependency.

\noindent First, following Leroux {\it et al.} \cite{GL96,LL01} we use  piecewise constant functions to 
appro\-xi\-mate  $b$ $(b'(X)=\sin\theta(X))$ as well as the term  $\cos\theta$ and the cross section area $S$.
Adding the equations $\partial_t Z=0$, $\partial_t \cos\theta=0$ and $\partial_t S=0$, the \textbf{PFS}-model can be written under a non-conservative form
with the variable $\WW =(b,\cos\theta,S,A,Q)^t$:
\begin{equation}\label{PFSncv}
\partial_t \WW+ \partial_X \FF(X,\WW) + B(X,\WW)\partial_X \WW =TS(\WW)
\end{equation}
where  $$\FF(X,\WW)=
\left(
\begin{array}{c}
0\\0\\0\\Q\\
\dsp \frac{Q^2}{A} + p(X,A)\end{array}
\right), \textrm{  }TS(\WW)=
\left(
\begin{array}{c}
0\\0\\0\\0\\
\dsp -g\,K(X,\SE)\, \frac{Q |Q|}{A}
\end{array}
\right)$$ and
$$B(X,\WW)=
\left(
\begin{array}{ccccc}
 0&0&0&0&0\\
 0&0&0&0&0\\
 0&0&0&0&0\\
 0&0&0&0&0\\
 gA&gA\overline{Z}&-c^2(A/\SE-1)- \mathcal{I} (X,\WW) &0&0
\end{array}
\right) \quad
$$
where we have written the pressure source term due to the geometry $ g I_2(X,\SE) \cos(\theta) $ as $\mathcal{I }(X,\WW) S'$. For instance, for a  circular cross-section pipe we have:
 $$
 \mathcal{I}(X,\WW) =\frac{1}{2\,\pi }\left(\dsp\frac{\mathcal{H}(\SE)\pi}{2}+
\mathcal{H}(\SE)\arcsin\left(\frac{\mathcal{H}(\SE)}{R(X)}\right)+\frac{\sigma(X,\mathcal{H}(\SE))}{2}\right)\,.
$$
\noindent Let $W_i^n$ be an approximation of the mean value of $\WW$ on the mesh
$m_i$ at time $t_n$. Since the values of $b,\cos\theta,S$ are known, integrating the above equations over $]X_{i-1/2},X_{i+\frac{1}{2}}[\times [t_n,t_{n+1}[ $, we can write a Finite Volume scheme as follows:

\begin{equation}\label{FVscheme}
\begin{array}{lll}
\WW_i^{n+1}  & = &\WW_i^n  \dsp -\dsp \alpha_i \left(\FF(\WW^*_{i+1/2}(0^-,\WW_i^n,\WW_{i+1}^n)) -
\FF(\WW^*_{i-1/2}(0^+,\WW_{i-1}^n,\WW_{i}^n))\right) \\
& & + \dsp TS(\WW^n_i)
\end{array}
\end{equation}

with $\alpha_i =\dsp \frac{\Delta t^n}{h_i} .$

\noindent  $\WW^*_{i+1/2}(\xi=x/t,\WW_i,\WW_{i+1})$ is  the exact or an approximate solution to the Riemann problem
at interface $X_{i+1/2}$ associated to the left and right states $\WW_i$ and $\WW_{i+1}$.  Let us also remark that the term $B(X,\WW)$ does not appear explicitly in this formulation since $b'$, $(\cos\theta)'$ and $S'$ are null on $]X_{i-1/2},X_{i+\frac{1}{2}}[$ but contributes to the computation of the numerical flux.

\noindent The computation of the interface quantities $\WW^*_{i\pm1/2}(0^{\pm},\WW_i,\WW_{i+1})$ will depend on two types of interfaces located at the point $X_{i+\frac{1}{2}}$ : the first one is a non transition point,  when the flow on both sides of the interface is of the same type. The second one is a transition point,
 when the flow changes of type through the interface. We recall the approach used in \cite{BG07} and adapt it here to the new terms.
According to the type of interface, we have to solve two different linearised Riemann problems.

\subsubsection{The Case of a non transition point}\label{SubSubSectionCaseOfANonTransitionPoint}
Expanding the term $\partial_X \FF(X,\WW)$ in the non-conservative equations (\ref{PFSncv}), the quasilinear formulation of the \textbf{PFS}-model (\ref{PFS}) reads:
\begin{equation}\label{FNC}
\partial_t \WW + D(\WW) \;\partial_X \WW =TS(\WW)
\end{equation} with $D$ the convection matrix defined by
\begin{equation}\label{ConvectionMatrix}
D(\WW) = \left(
\begin{array}{ccccc}
0 & 0 & 0 & 0 &  0\\
0 & 0 & 0 & 0 &  0\\
0 & 0 & 0 & 0 &  0\\
0 & 0 & 0 & 0 &  1\\
gA & gA\mathcal{H}(\SE) & \Psi(\WW) & c^2(\WW)-u^2 &  2u\\
\end{array}
\right)
\end{equation}

\noindent where $\dsp\Psi(\WW) = gS \partial_S \mathcal{H}(\SE) \cos\theta -c^2(\WW)\dsp\frac{A}{\SE}$ and $u=Q/A$ denotes the speed of the water. $c(\WW)$ is then  $c$ for the  pressurized flow  or  $\dsp \sqrt{g \frac{A}{T(A)}\cos\theta}$ for the free surface flow.

\begin{rque}

\noindent Let us remark that, since $\partial_X I_1(X,A) = I_2(X,A)+\partial_A I_1(A) \partial_X A$, the pressure source term $I_2$ does not appear in Equation (\ref{FNC}).
\end{rque}

\noindent To compute the interface quantities denoted by $(AM,QM)$ for the left hand side and
$(AP,QP)$ for the right hand side (see Figure \ref{QAMP} below), we solve  the following linearised Riemann problem:
\begin{equation}\label{rielin}
\left\lbrace
\begin{array}{lcl}
\partial_t \WW + \widetilde{D} \;\partial_X \WW & = & 0\\
\WW & = &
\left\lbrace
\begin{array}{lcr}
\WW_l=(b_l,\cos\theta_l,S_l,A_l,Q_l)^t & \hbox{ if }& x<0\\
\WW_r=(b_r,\cos\theta_r,S_r,A_r,Q_r)^t & \hbox{ if }& x>0
\end{array}
\right.
\end{array}
\right.
\end{equation}
\noindent with $(\WW_l,\WW_r)=(\WW_i,\WW_{i+1})$ and $\widetilde{D}=\widetilde{D} (\WW_l,\WW_r)=
\dsp D(\widetilde\WW) $ where $\dsp \widetilde{\WW}$ is some approximate state of the left $\WW_l$ and the right $\WW_r$  state.
\begin{rque}\label{RemarkChoiceOfTildeD}

\noindent  We will see in Section \ref{SectionWB} that the classical approximation $D(\widetilde{\WW})$ of the Roe matrix $\dsp
D_{Roe}(\WW_l,\WW_r) = \int_0^1 D(\WW_r+(1-s)(\WW_l-\WW_r))\,ds$ defined by $\widetilde{D} = \dsp D(\widetilde{\WW}) =
D\left(\frac{\WW_l+\WW_r}{2}\right)$ is not a suitable  choice to preserve the still water steady state.  However, we  propose in Section
\ref{SectionWB} a new approximation of $\widetilde{D}$ which  maintains it perfectly.
\end{rque}
\noindent We have then $W^*(0+,\WW_l,\WW_r) = (b_r,\cos\theta_r,S_r,AP,QP)^t$.

\noindent The eigenvalues of the matrix $\widetilde{D}$ are $\lambda_1=  0$, $\lambda_2=  0$, $\lambda_3=  0$, $\lambda_4 =\ut-c(\widetilde\WW) $, $\lambda_5 = \ut + c(\widetilde\WW)$ and the associated right eigenvectors:
$$
r_1(\widetilde\WW)=
\left(
\begin{array}{c}
c^2(\widetilde\WW)-\ut^2\\0\\0\\-g\At\\0
\end{array}
\right), \quad
r_2(\widetilde\WW)=
\left(
\begin{array}{c}
\Psi(\widetilde\WW)\\0\\-g\At\\0\\0
\end{array}
\right), \quad
r_3(\widetilde\WW)=
\left(
\begin{array}{c}
\mathcal{H}(\widetilde S)\\-1\\0\\0\\0
\end{array}
\right), \quad $$
$$
r_4(\widetilde\WW)=
\left(
\begin{array}{c}
0\\0\\0\\1\\\ut-c(\widetilde\WW)
\end{array}
\right), \quad
r_5(\widetilde\WW)=
\left(
\begin{array}{c}
0\\0\\0\\1\\\ut+c(\widetilde\WW)
\end{array}
\right).
$$
%
\begin{figure}[H]
\centering
\includegraphics[scale=0.5]{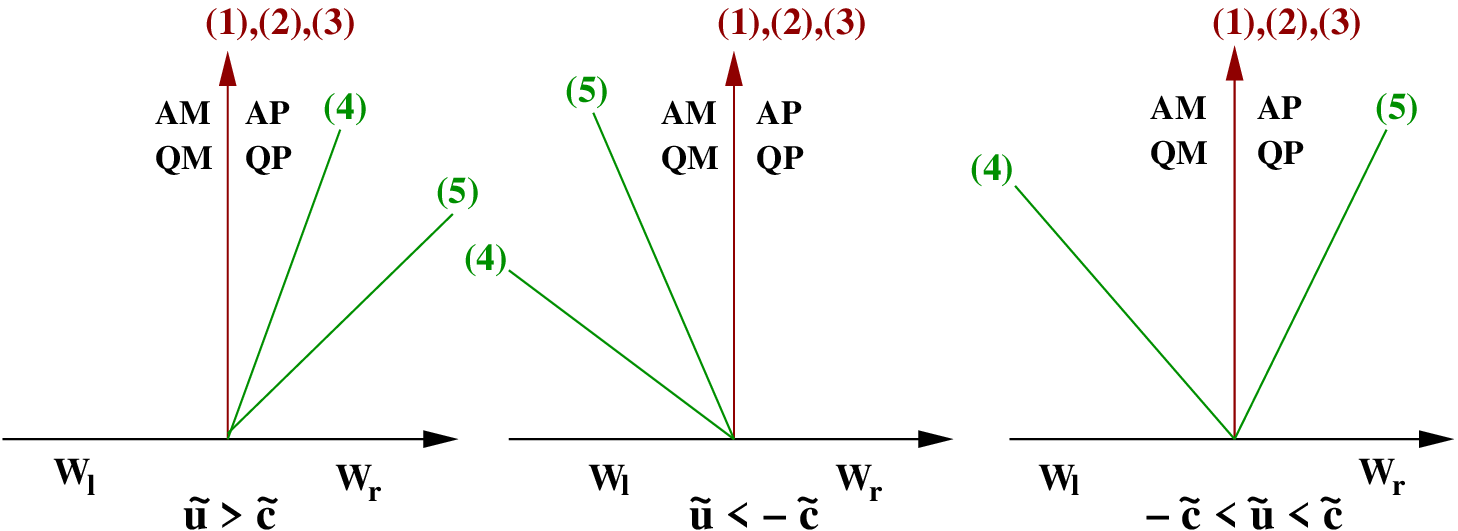}
\caption{Solution of the Riemann problem (\ref{rielin}). The lines $(i)$  corresponds
to the characteristic lines  $X/t=\lambda_i$, for $i=1,\ldots,5$ .}
\label{QAMP}
\end{figure}
\noindent We denote $P$ the transition matrix associated to the right eigenvectors of $\widetilde D$
and $P^{-1}$ its inverse. Setting $\jump{\WW} = \WW_r-\WW_l$, the solution of the Riemann problem are constant states connected by shocks propagating
along the characteristic lines $X/t=\lambda_i$. The jump associated to the eigenvectors $r_i$ is then equal to
$(P^{-1}\jump{\WW})_i\,r_i$. In particular, the discharge is continuous through the line $X/t=0$ since the fifth component of vectors
$r_1$, $r_2$ and $r_3$ are null. It follows that the equation associated to the wet area $A$ remains conservative.

\noindent Thus,  for instance in the subcritical case (when $-\ct<\ut<\ct$), we have:
\begin{equation*}
\begin{array}{l}
\dsp AM  =  A_l + \frac{g\,\At}{2\,\ct\,(\ct - \ut)}\,\psi_l^r + \frac{\ut + \ct}{2\,\ct}\,(A_r - A_l)
- \frac{1}{2\,\ct}\,(Q_r - Q_l)\\
      \\
\dsp QM=QP  =  Q_l - \frac{g\,\At}{2\,\ct}\,\psi_l^r + \frac{\ut^2 - \ct^2}{2\,\ct}\,(A_r - A_l)
- \frac{\ut - \ct}{2\,\ct}\,(Q_r - Q_l) \\
\\
\dsp AP  =  AM + \frac{g\,\At}{\ut^2 - \ct^2}\,\psi_l^r
\end{array}
\end{equation*}
where $\psi_l^r $ is the upwinded source term $
b_r-b_l+\mathcal{H}(\widetilde{\SE})(\cos\theta_r-\cos\theta_l)+\Psi(\widetilde{\WW})(S_r-S_l)$.
\begin{rque}

\noindent The friction term can also be upwinded in the same  way. Writing  the friction term under a conservative form
$$\partial_X \int_{X_0}^X \dsp \,K(s,\SE)\, \frac{Q(t,s) |Q(t,s)|}{A^2(t,s)}\, ds$$ (for some arbitrary $X_0$)
allows us to write the ``static'' slope $b$ as a ``dynamic'' one  as follows: $$b+\int_X \dsp \,K(s,\SE)\,
\frac{Q(t,s) |Q(t,s)|}{A^2(t,s)}\, ds$$ that we denote again $b$. Thus, the upwinding of the dynamic slope $b_{i+1}-b_i$ reads:
$$b_{i+1}-b_{i}+\dsp\int_{X_i}^{X_{i+1/2}} \frac{1}{K_s^2}\left\{ \frac{Q \abs{Q}}{A^2 R_h(\SE)^{4/3}}\right\}\,ds
+\int_{X_{i+1/2}}^{X_{i+1}} \frac{1}{K_s^2}\left\{ \frac{Q \abs{Q}}{A^2 R_h(\SE)^{4/3}}\right\}\,ds$$ which is equal to:
$$
b_{i+1}-b_{i} \dsp+ (X_{i+1/2}-X_i)  \frac{Q_{i} \abs{Q_{i}}}{ K_s^2\,A_{i}^2 R_h(\SE_{i})^{4/3}} + (X_{i+1}-X_{i+1/2}) \frac{Q_{i+1}
\abs{Q_{i+1}}}{K_s^2\,A_{i+1}^2 R_h(\SE_{i+1})^{4/3}}\,
$$ since $A$ and $Q$ are constant on each cells.

\noindent The terminology ``dynamic'' and ``static'' slope is used since one  is $(t,x)$-dependent while the other is only $x$-dependent.
\end{rque}

\subsubsection{Case of transition point}
In the case of a transition point, we assume that the propagation of the interface (pressurized-free surface or  free surface-pressurized) has a constant speed $w$
during a time step. The half line  $x = w\,t$ is then the discontinuity line of $\widetilde{D} (W_l,W_r)$.\\
Let us now consider $\U^-=(A^-,Q^-)$ and  $\U^+=(A^+,Q^+)$ the (unknown) states res\-pec\-tively on the left 
and on the right hand side
of the line $x = w\,t$ with \linebreak
$\dsp w~=~\frac{Q^+~-~Q^-}{A^+~-A^-}$. Both states $\U_l$ and  $\U^-$ (resp. $\U_r$ and $\U^+$) correspond  to the same type of flow.
Thus it makes sense to define the averaged matrices in each zone as follows:
\begin{itemize}
\item for $x < w\,t$, we set $\widetilde{D}_l=\widetilde{D} (\WW_l,\WW_r) = D(\widetilde{\WW}_l)$ for some approximation $\widetilde{\WW}_l$ which connects the state $\WW_l$ and $\WW^-$ (see Remark \ref{RemarkChoiceOfTildeD}).

\item for  $x > w\,t$, we set  $\widetilde{D}_r=\widetilde{D}(\WW_l,\WW_r) = D(\widetilde{\WW}_r)$ for some approximation $\widetilde{\WW}_l$ which connects the state $\WW^+$ and $\WW_r$ (see Remark \ref{RemarkChoiceOfTildeD}).
\end{itemize}

\begin{figure}[H]
\centering
\includegraphics[scale=0.5]{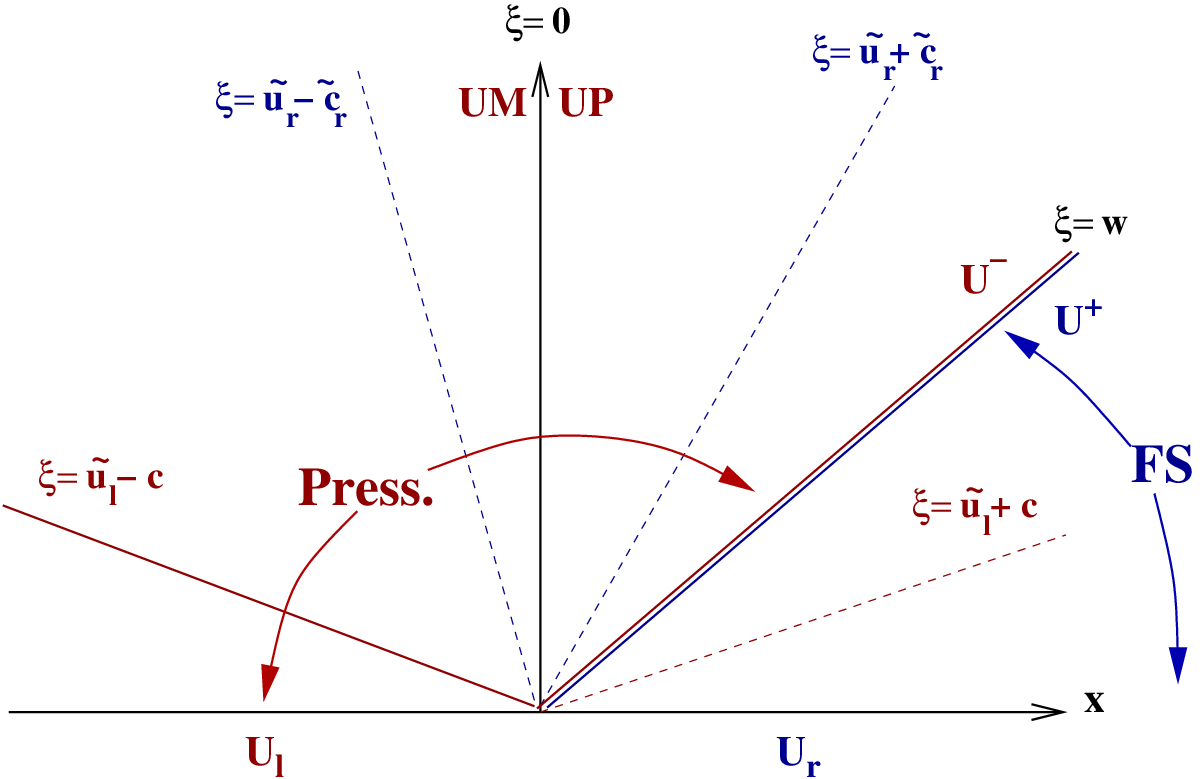}
\caption{Pressure state propagating downstream.}\label{plus1}
\end{figure}

\begin{figure}[!ht]
\centering
\includegraphics[scale=0.5]{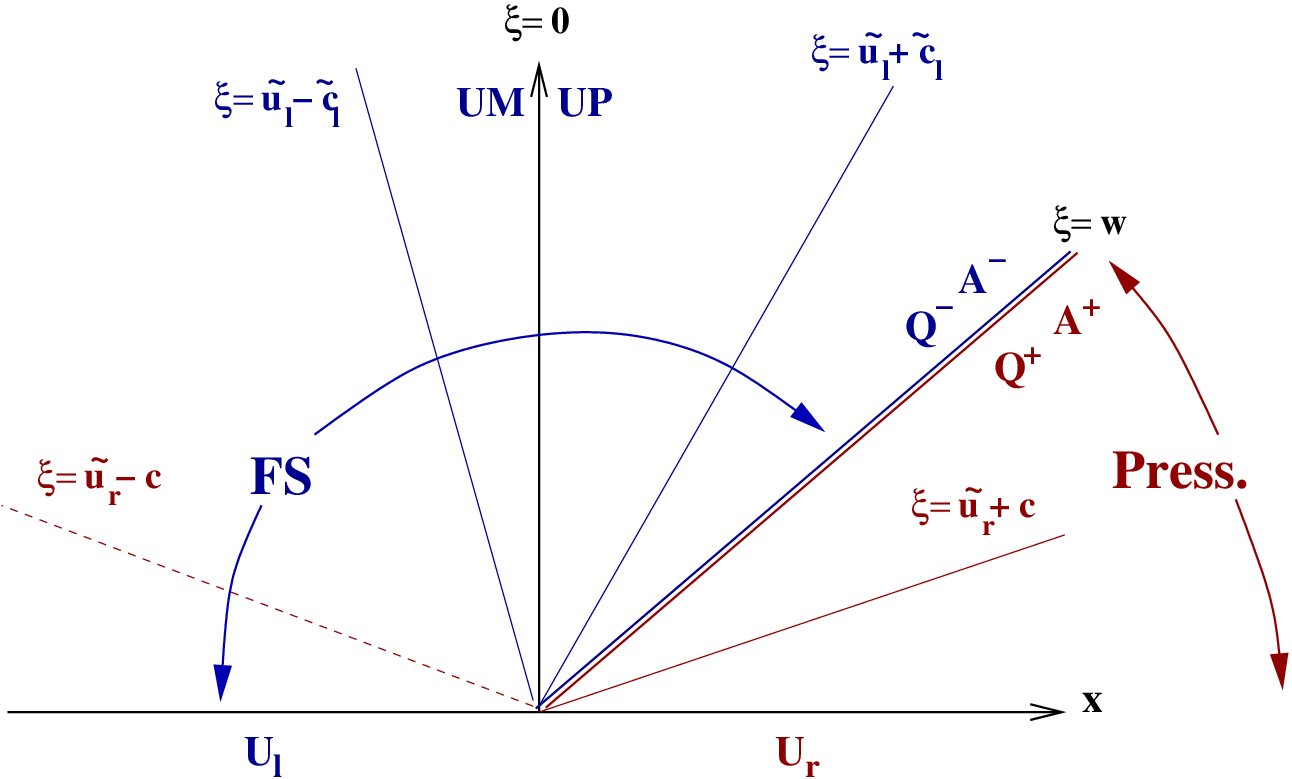}
\caption{Free surface state propagating downstream} \label{plus2}
\end{figure}

\noindent Then we formally solve two Riemann problems and use the Rankine-Hugoniot jump conditions through the line $x = w\,t$ which writes:
\begin{eqnarray}
Q^+ - Q^- & = & w\,(A^+ - A^-)\label{rh1}\\
F_5(A^+,Q^+) - F_5(A^-,Q^-)& =& w\,(Q^+ - Q^-)\label{rh2}
\end{eqnarray}
with $F_5(A,Q) = \dsp\frac{Q^2}{A} + p(X,A)$.
According to the left  ($\U^-$, $\UM$) and right unknowns ($\U^+$, $\UP$ ) at the interface $x_{i+1/2}$ and the sign of the speed $w$, we have to deal with four cases:
\begin{itemize}
 \item pressure state propagating downstream,
 \item pressure state propagating upstream,
 \item free surface state propagating downstream,
 \item free surface state propagating upstream.
\end{itemize}
We can next consider two  couples of ``twin cases'' : pressure state is propagating downstream (or upstream) as shown in the figure \ref{plus1} and free surface state propagating downstream (or upstream) as shown in the figure \ref{plus2}. Moreover, for all existing transition case, the upwinded altitude term $b_r-b_l$ in \cite{BG07} are replaced by $\psi_l^r$.

\paragraph{Pressure state propagating downstream:}
\noindent it is the case when on the left hand side of the line $\xi = w t$, we have a pressurized flow and on the right hand side we have a free surface flow: the speed $w$ of the transition point being positive.
Following Song \cite{SCL83} (see also \cite{M02}), an equivalent stationary hydraulic jump must occur from a supercritical to a subcritical condition and thus the characteristics speed satisfies the inequalities:
\begin{equation*}
\ut_r+\ct_r< w <\ut_l + c
\end{equation*}
where  $c$ is the sound speed for the pressure flow, $\ut_l$, $\ut_r$, and $\ct_r$ are defined by the same formula obtained
in the case of a non transition point but according to $\widetilde{D}_l$ and $\widetilde{D}_r$.

\noindent Therefore, only the characteristic lines drawn with solid lines are taken into account. Indeed they are related to incoming waves with respect to the cor\-responding space-time area $-\infty<\xi<w$. Conversely, the dotted line $ \xi=\ut_r-\ct_r$, for instance, related to the free surface zone but drawn in the area of pressurized flow is a ``ghost wave'' and is not considered. Thus $\U^+=\U_r$ and  $\U_l$, $\U^-$ are connected through the jumps across the characteristics $\xi=0$ and $\xi = \ut_l-c$.
Eliminating $w$ in the Rankine-Hugoniot jump relations (\ref{rh1})-(\ref{rh2}), we get $\U^-$ as the solution to the nonlinear system:

\begin{eqnarray}
\dsp (F_5(A_r,Q_r) - F_5(A^-,Q^-)) = \frac{(Q_r - Q^-)^2 }{(A_r - A^-)} \label{chslp1}\\
Q^- -Q_l - (A^- - A_l)(\ut_l - c) + \frac{g \psi^r_l \,\At_l}{c + \ut_l}  & = & 0\label{chslp2}
\end{eqnarray}
Finally, we obtain :
$$\left\{\begin{array}{lll}
 AP &=& A^- \\
 QM &=& Q^-\\
 QP &=& Q^-\\
 AM&=&AP-\dsp \frac{g\,\At_l\,\psi_l^r}{\ut_l^2 - c^2}.
 \end{array}\right.
$$

\paragraph{Free surface state propagating downstream:}
on the left hand side of the line $\xi = w t$ we have a free surface flow while on the right hand side, we have a pressurized flow (the speed $w$ of the transition point being positive). Following Song \cite{SCL83} again,  the characteristic speed satisfies the inequalities:
\begin{equation*}
\ut_l+\ct_l < w < \ut_r + c
\end{equation*}

\noindent There are two incoming characteristic lines with respect to the free surface area $-\infty<\xi<w$ (actually three with $\xi=0$) and they can connect the given left state $\U_l$ with any arbitrary free surface state $\UM$.  Thus only one  characteristic line ($\xi= \ut_r + c$) gives any information (it is  the equation (\ref{slchp1}) above) as an incoming characteristic line with respect to the pressurized zone $w<\xi<+\infty$. From the jump relations through the characteristic $\xi~=~0$, and after the elimination of $w$ in the Rankine-Hugoniot jump relations
(\ref{rh1}),(\ref{rh2}) we get another equation, namely Equation (\ref{slchp2}) above.
It remains to close the system of four unknowns $(A^-,\,Q^-,\,A^+,\,Q^+)$.
Firstly, we use a jump relation across the transition point (with speed $w$) for the total head $\dsp \Psi = \frac{u^2}{2} + c^2 \ln\left(\frac{A}{\SE}\right) + g\,\mathcal{H}(A)\cos\theta + g\,b$, from Equation (\ref{ThmPFSEquationForU}), which writes:
$$
\Psi^+ - \Psi^- = w\,(u^+ - u^-)\,.
$$
Finally, we use the relation:
$$ w = w_{pred} \, \textrm{ with } w_{pred} = \frac{Q_r-Q_l}{A_r-A_l}\,.$$
We have then to solve the nonlinear system:
\begin{eqnarray}
(Q_r - Q^+)  = (A_r - A^+)\,(\ut_r + c) \label{slchp1}\\
\nonumber\\
(Q^+ - Q^-)\,(Q_r - Q_l)  =  (A_r - A_l)\, (F_2(A^+,Q^+) - F_2(A^-,Q^-))\label{slchp2}\\
\nonumber\\
\frac{(Q^+)^2}{2\,(A^+)^2} + c^2 \ln\left(A^+\right) + g\cos\theta\,\mathcal{H}(A^+)
- \frac{(Q^-)^2}{2\,(A^-)^2}- c^2 \ln\left(A^-\right)- g\cos\theta\,\mathcal{H}(A^-) \nonumber\\
=  \frac{Q_r - Q_l}{A_r - A_l}\,\left( \frac{Q^+}{A^+} - \frac{Q^-}{A^-}  \right) \label{slchp3}\\
\nonumber\\
(Q_r - Q_l)\, (A^+ - A^-) = (Q^+ - Q^-)\, (A_r - A_l)\label{slchp4} %
\end{eqnarray}

\noindent The states $\UM$ et $\UP$ are then obtained by the following identities:
\begin{equation*}
\begin{array}{l}
\dsp AM  =  A_l + \frac{g\,\At_l\,\psi_l^r}{2\,\ct_l(\ct_l - \ut_l)} +
 \frac{\ut_l+\ct_l}{2\,\ct_l}\,(A^- - A_l) - \frac{1}{2\,\ct_l}\,(Q^- - Q_l) \\[5mm]
\dsp AP =  AM + \frac{g\,\At_l\,\psi_l^r}{\ut_l^2-\ct_l^2}\\[5mm]
\dsp QM = QP = QMP = Q_l + \frac{g\,\At_l\,\psi_l^r}{2\,\ct_l} \;+ \\
\hspace{4.8cm}\dsp + \;\frac{\ut_l^2-\ct_l^2}{2\,\ct_l}\,(A^- - A_l) -
\frac{\ut_l - \ct_l}{2\,\ct_l}\,(Q^- - Q_l)
\end{array}
\end{equation*}
\noindent Finally, the update state $A_i^{n+1} \,,\, Q_i^{n+1}$ are obtained by the same relation as in the case of a non transition point.

\noindent Using equations (\ref{FVscheme}) we update the values of  $A_i^{n+1},\, Q_i^{n+1}$ with a standard stability condition of Courant-Friedrich-Levy  controlling the time step size $\Delta t^n$.

\subsubsection{Updating the state of the flow $\boldsymbol{E}$ in a  cell.}
\begin{figure}[H]
 \begin{center}
 \begin{tabular}{cc}
 \includegraphics[scale =0.25]{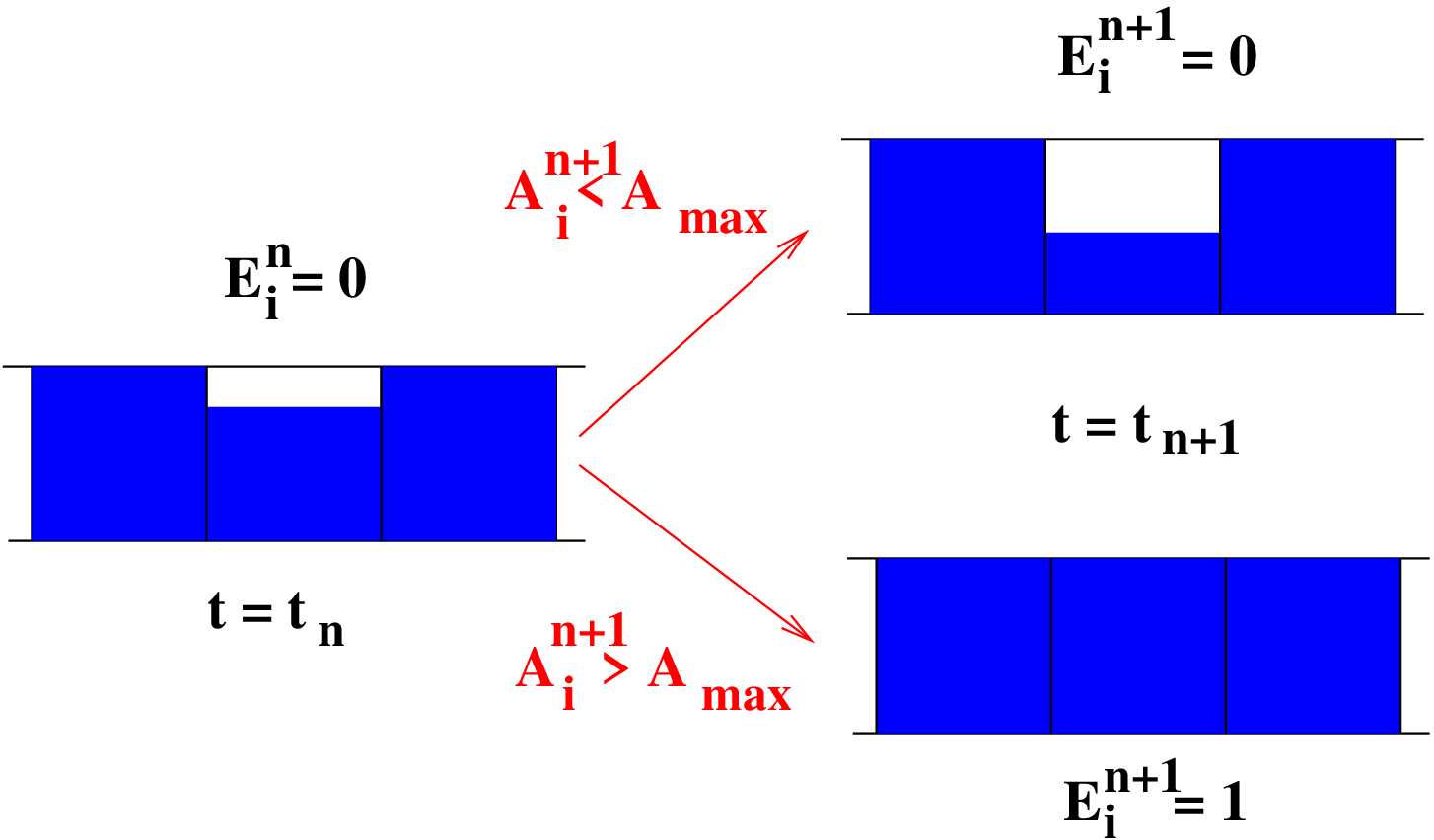} &
 \includegraphics[scale =0.25]{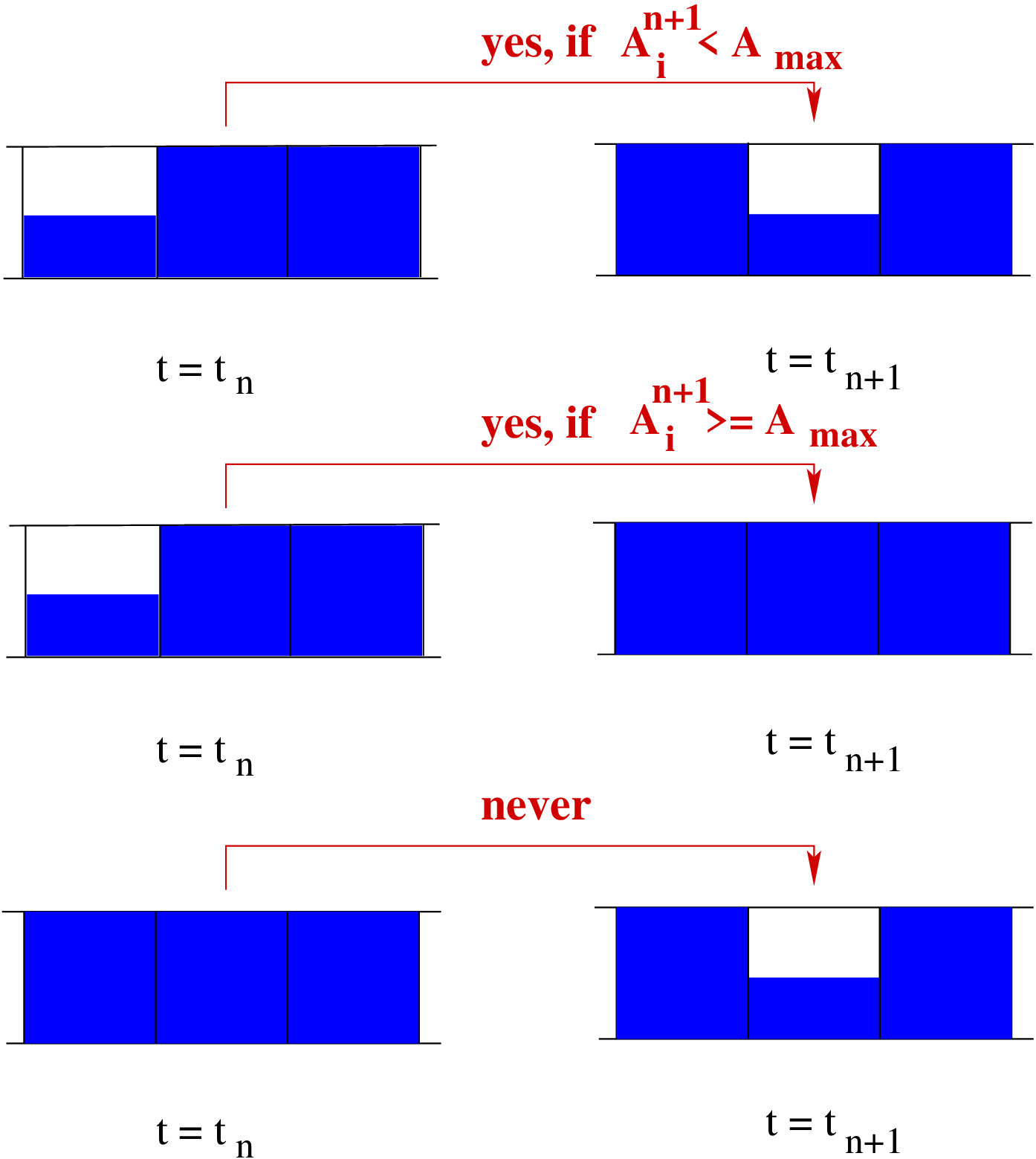}
 \end{tabular}
  \caption{Update of the state $E_i^{n+1}$ of the mesh $m_i$.}
  \label{UpdateA}
 \end{center}
\end{figure}

\noindent To update the state $E$ in the cell $m_i$ (see \fig\ref{UpdateA}), we use a discrete version of the state indicator $E$ defined by (\ref{DefinitionE}) equal to $1$ for a
pressurized flow and $0$ otherwise. Following \cite{BG07}, after the computation of the wet area $A_i^{n+1}$ we predict the state of the flow in the cell $m_i$ by the following criterion:
\begin{itemize}
\item[$\bullet$]  if $E_i^n = 0$ then :\\
if $A_i^{n+1}<S_i$ then   $E_i^{n+1}=0$, else $E_i^{n+1}=1$,

\item[$\bullet$]   if $E_i^n=1$ :\\
if   $A_i^{n+1}\geq S_i$ then   $E_i^{n+1}=1$,
else $E_i^n=E_{i-1}^n\cdot E_{i+1}^n$.
\end{itemize}

\noindent Indeed, if $A_i^{n+1}\geq S_i$ it is clear that the state of the flow in  the cell $m_i$ becomes pressurized, on the other hand if $A_i^{n+1}<{S_i}$ in a mesh previously pressurized, we do not know \textit{a priori} if the new state is free surface ($\rho=\rho_0$ and the value of the wetted area is less than ${S_i}$) or pressurized (in depression, with $\rho<\rho_0$ and the value of the wetted area is equal to $S_i$: see Remark \ref{RemarkDepressionOverpressureStates} and \fig\ref{RhoCriterion}).\\
So far as we do not take into account complex phenomena such that entrapment of air pockets or cavitation and keeping in mind that the  CFL condition  ensures that a transition point crosses at most one mesh at each time step, we  postulate that:
\begin{enumerate}
\item if the state of the flow in  the cell $m_i$ is free surface at time $t_n$, its state at time $t_{n+1}$ is only determined by the value of $A_i^{n+1}$ and it cannot become in depression.
\item if the state of the flow in  the cell $m_i$ is pressurized at time $t_n$ and if $A_i^{n+1}<S_i$, it becomes free surface if and only if at least one adjacent cell was free surface  at time $t_n$. This is exactly the discrete version of the continuous $\dsp\frac{A}{\SE}$ criterion in Remark \ref{RemarkDepressionOverpressureStates} and displayed on \fig\ref{RhoCriterion}.
\end{enumerate}

\begin{rque}\label{RemarkDepressionOverpressureStates}

\noindent As we do not take into account complex phenomena such that entrapment of air pockets,  each connected component of the pressurized 
area is simply connected (see  \fig\ref{RhoCriterion}). Moreover, for each depression area $D$, its closure $\overline{D}$ is a strict subset of the 
pressurized set. It follows that when $A<S$ on  each pressurized area,  we observe a depression as displayed on \fig\ref{RhoCriterion}. 
Moreover, we may also use a visual depression indicator given by the function $\dsp\frac{A}{\SE}$: the case $\SE=A$ corresponds to a free surface state while $\SE>S$ to an overpressure state and $\SE<S$ to a depression state. On  \fig\ref{RhoCriterion}, we draw the behavior of the interface speed $w$ in the $(X,t)$-plane and the graph of the function $\dsp\frac{A}{\SE}$ at fixed time $t_3$.
\begin{figure}[H]
 \begin{center}
 \includegraphics[scale =0.49]{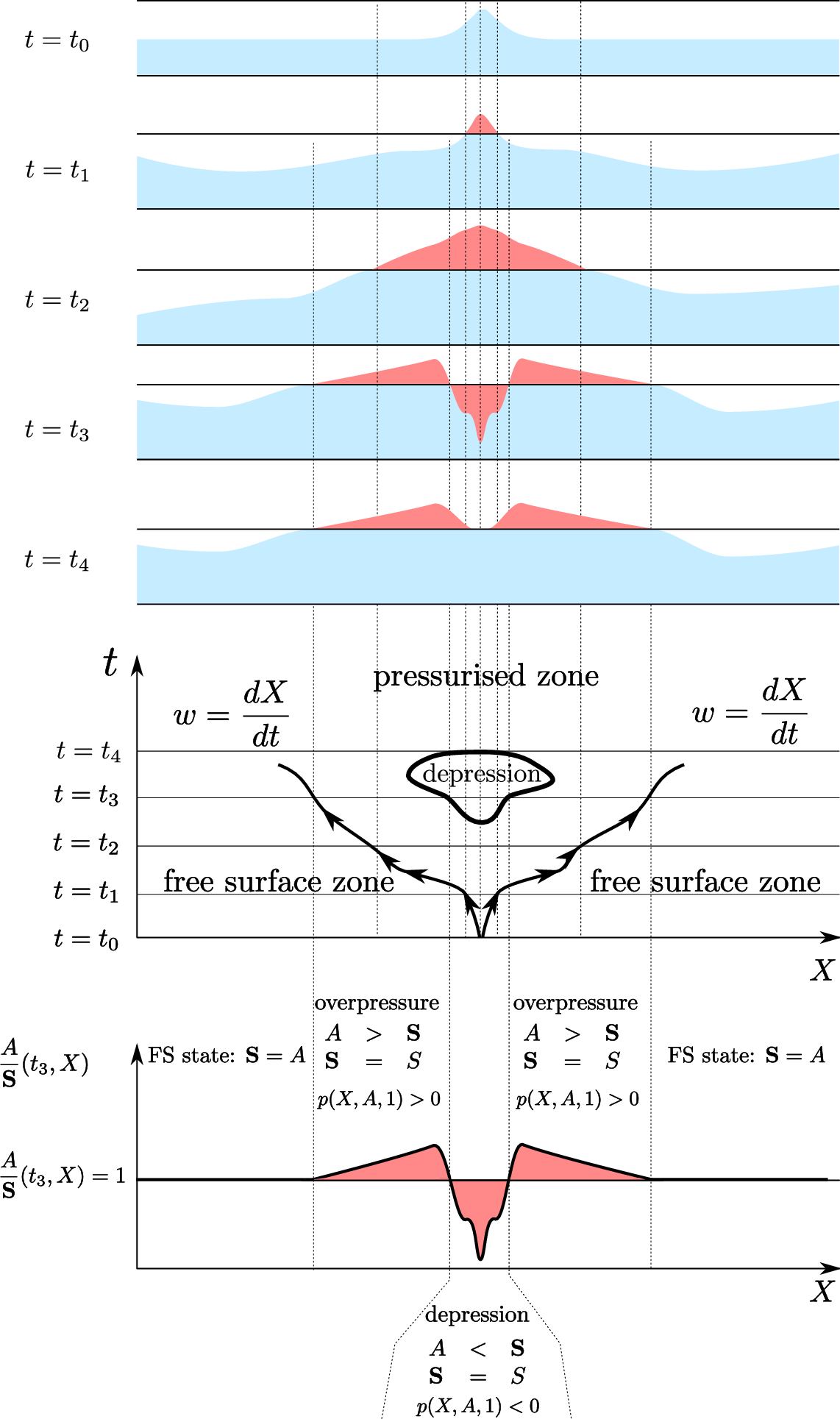}
  \caption{$\dsp\frac{A}{\SE}$ depression indicator.}
  \label{RhoCriterion}
 \end{center}
\end{figure}
\end{rque}

\section{Remarks on  still water steady state: an exactly well balanced scheme}\label{SectionWB}
This section is devoted to the construction of an exactly well-balanced scheme in the sense that it
maintains perfectly the still water steady states. This scheme, noted EWBS, is obtained by a suitable
definition of the convection matrix.

\noindent The numerical approximation of the \textbf{PFS}-model (\ref{PFSModelPressureInFlux}) reads:
\begin{eqnarray}
A_i^{n+1}  & = & A_i^n  \dsp -\frac{\Delta t}{h_i} \left(Q_{i+1/2}^n - Q_{i-1/2}^n\right)
\label{NumericalSchemeA} \\ \nonumber\\
Q_i^{n+1}  & = & Q_i^n  \dsp -\frac{\Delta t}{h_i} \left(F_5(AM_{i+1/2}^n,Q_{i+1/2}^n) -
F_5(AP_{i-1/2}^n,Q_{i-1/2}^n)\right) \label{NumericalSchemeQ}
\end{eqnarray} where $Q_{i\pm 1/2}$ stands for $QMP_{i\pm 1/2}$ and  $\dsp F_5(A,Q) = \frac{Q^2}{A} + c^2 (A-\SE) + g\,I_1(X,\SE)$.
For instance, in the subcritical case, the interface quantities reads:
\begin{equation}\label{InterfaceQuantitiesSubcriticalCase}
\begin{array}{lll}
\dsp AM_{i+1/2}^n  &=& \dsp A_i^n + \frac{g\,\At_{i+1/2}^n}{2\,\ctt_{i+1/2}\,(\ctt_{i+1/2} - \ut_{i+1/2}^n)}\,\psi_i^{i+1} \\ \\& &\dsp +
\frac{\ut_{i+1/2}^n + \ctt_{i+1/2}}{2\,\ctt_{i+1/2}}\,(A_{i+1}^n - A_i^n)  -
 \frac{1}{2\,\ctt_{i+1/2}}\,(Q_{i+1}^n - Q_i^n)\\
      \\
\dsp Q_{i+1/2}^n  &=&  \dsp Q_i^n - \frac{g\,\At_{i+1/2}^n}{2\,\ctt_{i+1/2}}\,\psi_i^{i+1} + \frac{\ut_{i+1/2}^{n,2} -
\ctt_{i+1/2}^2}{2\,\ctt_{i+1/2}}\,(A_{i+1}^n - A_i^n) \\ \\ & & \dsp- \frac{\ut_{i+1/2} - \ctt_{i+1/2}}{2\,\ctt_{i+1/2}}\,(Q_{i+1}^n -
Q_i^n)
\\
\\
\dsp AP_{i+1/2}^n & =& \dsp AM_{i+1/2}^n + \frac{g\,\At_{i+1/2}^n}{ \ut_{i+1/2}^{n,2} - \ctt_{i+1/2}^2}\,\psi_i^{i+1}
\end{array}
\end{equation}  where the upwinded source term reads:
$$b_{i+1}-b_i+\mathcal{H}(\widetilde{\SE}_{i+1/2}^n)(\cos\theta_{i+1}-\cos\theta_i)+\Psi(\widetilde{\WW}{i+1/2}^n)(S_{i+1}-S_i)$$ and
$\ctt_{i+1/2}$ stands for $ c(\widetilde{\WW}_{i+1/2}^n)$ with $$\widetilde{\WW}_{i+1/2}^n = \left(\widetilde{b}_{i+1/2},\,
\widetilde{\cos\theta}_{i+1/2},\, \widetilde{S}_{i+1/2},\, \At_{i+1/2}^n,\, \Qt_{i+1/2}^n\right)$$  given by
\begin{equation}\label{DefinitionOfTildebThetaSQ}
\begin{array}{l}
\dsp \widetilde{b} = \frac{b_i+b_{i+1}}{2},\, \dsp \widetilde{\cos\theta} = \frac{\cos\theta_i+\cos\theta_{i+1}}{2},\,
 \dsp \widetilde{S} = \frac{S_i+S_{i+1}}{2},\, \Qt_{i+1/2}^n = \frac{Q_i^n+Q_{i+1}^n}{2},
\end{array}
\end{equation} and the approximation of $\At_{i+1/2}^n$ to be specified.

\noindent Starting from a discrete state $(A_i^n,Q_i^n)$ at time $t_n$ such that:
\begin{description}
 \item[(H)]let $n$ such that: $\forall i$, $Q_i^n = 0$ and $A_i^n$ satisfy the discrete still
water steady state  equation (according to Equation (\ref{ThmPFSEquationForU})):
\begin{equation}\label{DiscreteSteadyStateEquation}
c^2\ln\left(\frac{A_{i+1}^n}{S_{i+1}}\right)+g \mathcal{H}(\SE_{i+1}^n)\cos\theta_{i+1} +g b_{i+1}
=c^2\ln\left(\frac{A_{i}^n}{S_{i}}\right)+g \mathcal{H}(\SE_{i}^n)
\cos\theta_i +gb_{i}\,,
\end{equation}
\end{description}
we will say that:
\begin{defi}\label{DefinitionWBS}
\begin{enumerate}
\item[]
\item  The numerical scheme (\ref{NumericalSchemeA}-\ref{NumericalSchemeQ})-(\ref{DefinitionOfTildebThetaSQ}) for some approximations of the
terms $\At_{i\pm 1/2}^n$ is
$(k_A,k_Q)$ well-balanced (also denoted by  $(k_A,k_Q)$-WB) if:
$$
\forall i,\,\begin{array}{l}
\abs{A_i^{n+1}-A_i^{n}} = O((\max_{i\in\Z} h_i)^{k_{A}}) \textrm{ and } \abs{Q_i^{n+1}-Q_i^{n}} = O((\max_{i\in\Z} h_i)^{k_{Q}})\,,
\end{array}
$$
with $k_A>1$, $k_Q>1$ the well-balanced order of the numerical scheme (\ref{NumericalSchemeA}) and (\ref{NumericalSchemeQ}) respectively.
\item  The numerical scheme (\ref{NumericalSchemeA}-\ref{NumericalSchemeQ})-(\ref{DefinitionOfTildebThetaSQ}) for some approximations of the
terms $\At_{i\pm 1/2}^n$ is exactly
well-balanced (also denoted by EWB) if:
$$
\forall i,\,\begin{array}{l}
\abs{A_i^{n+1}-A_i^{n}} =0 \textrm{ and } \abs{Q_i^{n+1}-Q_i^{n}} =0\,.
\end{array}
$$
\end{enumerate}
\end{defi}
We will denote by $(k_A,k_Q)$-WBS the $(k_A,k_Q)$ well-balanced scheme and EWBS the exactly well-balanced scheme.

\begin{description}
\item[(SF)]  In the rest of this paper, we assume $h_i = \Delta X$ constant, the radius
$R$ and $b$ are given  affine functions, the angle
$\theta$  is constant which implies that the jumps across the interface $X_{i+1/2}$: $\Delta R_{i+1/2} = R_{i+1}-R_{i}=\Delta R$, $\Delta b_{i+1/2}
=b_{i+1}-b_{i}=\Delta b$ are constant and $\Delta\cos\theta_{i+1/2}$ is null.
\end{description}
Then under this simplified framework, we show that:
\begin{thm}\label{ThmWB}
\begin{enumerate}
 \item[]
\item
\begin{itemize}
 \item The numerical scheme (\ref{NumericalSchemeA}-\ref{NumericalSchemeQ})-(\ref{DefinitionOfTildebThetaSQ}) with the classical choice
 \begin{equation}\label{DemiSomme}
 \At_{i+1/2}^n =\dsp \frac{A_i^n+A_{i+1}^n}{2}
 \end{equation} and non constant section $S$ is not  well-balanced in the sense of Definition \ref{DefinitionWBS} (we have $k_Q=1$).
 \item For constant section and $Z=0$, the numerical scheme
(\ref{NumericalSchemeA}-\ref{NumericalSchemeQ})-(\ref{DefinitionOfTildebThetaSQ}) with  (\ref{DemiSomme}) is EWB.
 \end{itemize}
\item Under a suitable choice of $\At_{i\pm1/2}^n$, the numerical
scheme (\ref{NumericalSchemeA}-\ref{NumericalSchemeQ})-(\ref{DefinitionOfTildebThetaSQ}) is EWB.
\end{enumerate}
\end{thm}
The following section deals with the numerical scheme with
$\At_{i+1/2}^n$ defined as (\ref{DemiSomme}) where we show the first point of  Theorem \ref{ThmWB}. The second point  is studied in
Section \ref{SubSectionEWBS} where a convenient definition of $\At_{i+1/2}^n$ leads to an EWBS.

\subsection{Still water steady state and the classical approximation}\label{SubSectionSteadyStatesWithClassicalApproximation}
The simpler choice of definition for the convection matrix $D(\widetilde{\WW})$ is the one obtained by the approximation
of the mean value of the Roe matrix $\dsp D_{Roe}(\WW_l,\WW_r) = \int_0^1 D(\WW_r+(1-s)(\WW_l-\WW_r))\,ds$. This approximation is given by
$\widetilde{D} = \dsp
D(\widetilde{\WW}) = D\left(\frac{\WW_l+\WW_r}{2}\right)$ that we call ``the classical approximation''. Thus, defining $\At_{i+1/2}^n$
as follows: $$\At_{i+1/2}^n =\dsp
\frac{A_i^n+A_{i+1}^n}{2}\,$$ provides the classical approximation.  But, it is not
suitable to preserve the still water steady state: we will see that the
numerical scheme (\ref{NumericalSchemeA}-\ref{NumericalSchemeQ})-(\ref{DefinitionOfTildebThetaSQ}) with (\ref{DemiSomme}) defines a non
well-balanced scheme in the sense of  Definition (\ref{DefinitionWBS}) since $k_Q=1$.

\noindent To this end, let us assume \textbf{(H)} and  \textbf{(SF)} at time $t_n$, 
Equations (\ref{InterfaceQuantitiesSubcriticalCase}) read for every $i$:
\begin{equation}\label{AMAPQMP}
\begin{array}{ll}
\displaystyle AM_{i+1/2}^n  &= \dsp A_i^n + \frac{g\,\At_{i+1/2}^n}{2\,\ctt_{i+1/2}^2}\,\psi_i^{i+1} +
\frac{\Delta A_{i+1/2}^n}{2}\\ \\
\displaystyle Q_{i+1/2}^n  &=\dsp  -\frac{g\,\At_{i+1/2}^n}{2\,\ctt_{i+1/2}}\,\psi_i^{i+1} -
\ctt_{i+1/2}\frac{\Delta A_{i+1/2}^n}{2} \\ \\
\displaystyle AP_{i+1/2}^n  &= \dsp AM_{i+1/2}^n -\frac{g\,\At_{i+1/2}^n}{\ctt_{i+1/2}^2}\,\psi_i^{i+1}
\end{array}\,.
\end{equation}

\noindent Denoting $Q_{i+1/2}-Q_{i-1/2}$ by $\Delta Q_{i+1/2}$, we have:
$$
\begin{array}{lll}
\Delta Q_{i+1/2}^n &=& \dsp \frac{g}{2\, \ctt_{i-1/2}\, \ctt_{i+1/2} }
\Big\{
\left(\psi_{i-1}^{i}-\psi_{i}^{i+1} \right) \ctt_{i+1/2}\,\At_{i-1/2}^n
\\ & & \\
& & \dsp - \left(\At_{i+1/2}^n - \At_{i-1/2}^n\right)
\ctt_{i+1/2}\,\psi_{i}^{i+1}
\\ & & \\
& & \dsp + \left(\ctt_{i+1/2} -
\ctt_{i-1/2}\right)\At_{i+1/2}^n\,\psi_{i}^{i+1}
\Big\}
\\ & & \\
& & \dsp  + \frac{\Delta A_{i+1/2}^n}{2} \left( \ctt_{i-1/2} -
\ctt_{i+1/2}\right)
\end{array}
$$
\noindent where
$$
\begin{array}{lll}
\dsp \At_{i+1/2}^n-\At_{i-1/2}^n &=&\dsp   \Delta A_{i+1/2}^n \\ & & \\
\dsp  \psi_{i-1}^{i}-\psi_{i}^{i+1} &=&
\left\{
\begin{array}{l}
\dsp-\Delta b \quad\textrm{ if }  E_{i}= 0 \\ \\ \\
\dsp-\Delta b -\frac{g\,\cos\theta}{2}\Delta R  \\ -\dsp\frac{c^2\,\Delta S_{i+1/2}}{S_{i-1/2}\,S_{i+1/2}} \\ \times \left( \Delta S_{i+1/2}\,
\At_{i-1/2}^n 
-\Delta A_{i+1/2}^n\,S_{i-1/2}\right) \textrm{ if }  E_{i} = 1
\end{array}
\right.\\ & & \\
\dsp \abs{\ctt_{i+1/2} -\ctt_{i-1/2}} & &
\left\{
\begin{array}{ll}
\leqslant & C \Delta X  \textrm{ if } E_{i} = 0\,\textrm{ ( for some constant } C) \\
=  & 0 \textrm{ if } E_{i} = 1\,(\textrm{ since } \ctt_{i+1/2}=\ctt_{i-1/2} = c)
\end{array}\right.
\end{array}.
$$
Denoting then $$M = \max\left( \max_i{\left(\ctt_{i+1/2}\,\At_{i-1/2}^n\right)} ,
\max_i{\left(\ctt_{i+1/2}\,\abs{\psi_{i}^{i+1}}\right)} ,
\max_i{\left(\At_{i+1/2}^n\,\abs{\psi_{i}^{i+1}}\right)} , C \right)$$ and observing that $$\forall i,\,O(\Delta R)=O(\Delta
S_{i+1/2})=O(\Delta A_{i+1/2}^n)=O(\Delta X)$$ and $$O(\ctt_{i+1/2}\,\ctt_{i-1/2})=O(S_{i+1/2}\,S_{i-1/2})=O(1)\,,$$ we deduce:
$$\abs{\Delta Q_{i+1/2}^n} \leqslant M {\Delta  x^2}\,.$$ It follows that:
$$\abs{A_{i}^{n+1}-A_{i}^n} = O({\Delta X^2})\,.$$

\noindent Then, we denote $F_{5}(AM_{i+1/2}^n,Q_{i+1/2}^n)-F_{5}(AP_{i-1/2}^n,Q_{i-1/2}^n)$ by $\Delta F = T_1+T_2$
with

$$T_1
=\dsp\frac{\left(AM_{i+1/2}^n-AP_{i-1/2}^n\right)\left(Q_{i+1/2}^n\right)^2+AM_{i+1/2}^n\left(\left(Q_{i+1/2}^n\right)^2-\left(Q_{i-1/2}
^n\right)^2\right)}{AM_{ i+1/2}^n \,AP_{i-1/2}^n}$$ and
$$
T_2 =
\left\{
\begin{array}{lll}
 \dsp g\cos\theta\left(I_1(X_{i+1/2,},AM_{i+1/2}^n)-I_1(X_{i+1/2,},AP_{i-1/2}^n)\right) & \textrm{ if } & E_{i} = 0\\
& &\\
c^2(AM_{i+1/2}^n-AP_{i-1/2}^n) + c^2\Delta S_{i+1/2}                        & \textrm{ if } & E_{i} =1
\end{array}
\right.
$$
\noindent where $\dsp\abs{\left(Q_{i+1/2}^n\right)^2-\left(Q_{i-1/2}
^n\right)^2} =O(\Delta X^3)$ and
\begin{equation}\label{DiffAInterface}
\begin{array}{lll}
 {AM_{i+1/2}^n-AP_{i-1/2}^n} &=& \dsp \Delta A_{i+1/2}^n +
\frac{g}{2}\left(\frac{\At_{i+1/2}^n\psi_i^{i+1}}{\ctt_{i+1/2}^2}+\frac{\At_{i-1/2}^n\psi_{i-1}^{i}}{\ctt_{
i-1/2 }^2}\right)
\end{array}\,.
\end{equation}
As the term $\dsp\left(\frac{\At_{i+1/2}^n\psi_i^{i+1}}{\ctt_{i+1/2}^2}+\frac{\At_{i-1/2}^n\psi_{i-1}^{i}}{\ctt_{
i-1/2 }^2}\right)$ is at least of order $\dsp\Delta X$ since $$\psi_i^{i+1} = O(\Delta X),$$ we
have: $$\dsp\abs{AM_{i+1/2}^n-AP_{i-1/2}^n}=O(\Delta X)\,.$$
It follows that :$$\abs{Q_{i}^{n+1}-Q_{i}^n} =O({\Delta X})\,.$$

\begin{rque}
For constant section $S$ and $Z=0$, it is  easy to see that the numerical scheme
(\ref{NumericalSchemeA}-\ref{NumericalSchemeQ})-(\ref{DefinitionOfTildebThetaSQ}) with
(\ref{DemiSomme}) is  EWB.
\end{rque}

\noindent Although the scheme (\ref{NumericalSchemeA}-\ref{NumericalSchemeQ})-(\ref{DefinitionOfTildebThetaSQ}) with
(\ref{DemiSomme}) have an order $k_Q=1$ for non constant section, the still
water steady state for the pressurized case is very well maintained for great value of the sonic speed $c$. But it is
not the case for the free surface numerical scheme (see \fig\ref{StatioSL} with $\Delta X=10^{-3}$, $b=10^{-2}$, uniform pipe with diameter
$1$). We plot in \fig\ref{ComparaisonMeanMatrixStationnaryMatrixc30} and \fig\ref{ComparaisonMeanMatrixStationnaryMatrixc100} a still water
steady
state for two values of $c$ for a  given $\Delta X =1$ ( with $b=-0.9$ ): the absolute error obtained is
$10^{-5}$ for $c=30$ and $10^{-9}$ for
$c=200$. Indeed, let us consider, for the sake of simplicity, the case $S=cte$. At the continuous
level, the still water steady state equation reads: $$c^2\ln\left({A}\right)+g R\cos\theta +g b = cte.$$ With the hypothesis
\textbf{(H)}, in particular using  Equation (\ref{DiscreteSteadyStateEquation}), we write: $$A_{i+1}^n = A_i^n
\exp\left(-\frac{g}{c^2}\Delta b
\right)\,.$$
Given $\Delta X$, for great value of $c$, we can approximate $A_{i+1}^n$ by:
$$A_{i+1}^n \approx A_{i}^n\left(1 -\frac{g}{c^2}\Delta b\right)\,.$$
Then, replacing the right hand-side of $A_{i+1}^n$ in (\ref{DiffAInterface}), we have $$\abs{AM_{i+1/2}^n-  AP_{i-1/2}^n} = \,g\,A_i^n
\frac{\abs{\Delta b}}{c^2}=O\left(\frac{\Delta X}{c^2}\right)\,.$$ As we can see $k_Q$ is always equal to $1$ but the constant
$\dsp\frac{1}{c^2}$ plays the role of  a smoothing term  which helps the scheme
(\ref{NumericalSchemeA}-\ref{NumericalSchemeQ})-(\ref{DefinitionOfTildebThetaSQ}) with (\ref{DemiSomme}) to stabilize rapidly towards
the equilibrium.  Physically, $c$ is approximatively $1400$ (for a pressurized
flow
without air), thus $\dsp\frac{1}{c^2}\approx 1.9\,10^{-6}$. Since $c(A)\ll c$,  this feature is not observed for the free surface numerical
scheme.

\subsection{An exactly well-balanced scheme}\label{SubSectionEWBS}
\begin{figure}[H]
\begin{center}
\includegraphics[scale=0.7]{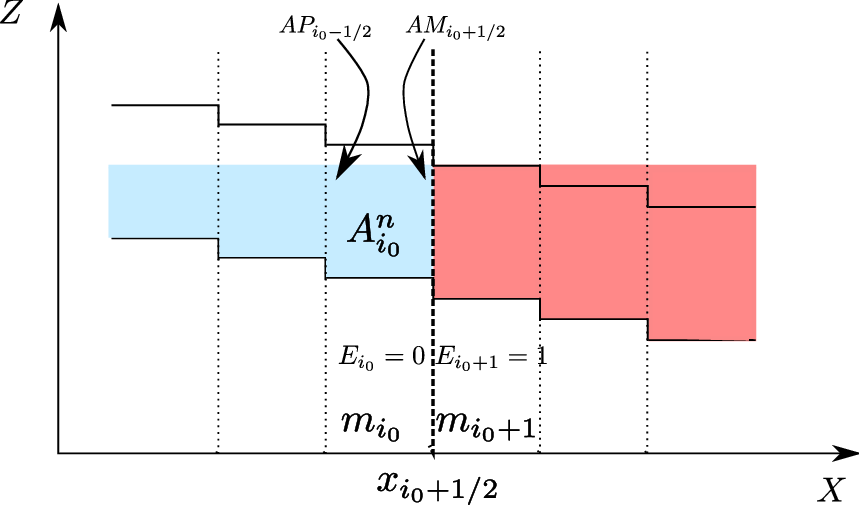}
\end{center}
\caption{Discrete representation of the mixed free surface-pressurized  still
water steady state at time $t^n$.}
\label{DiscretRepresentationDomain}
\end{figure}

This section is devoted to the construction of an EWBS. We have seen in Subsection \ref{SubSectionSteadyStatesWithClassicalApproximation}
that the classical approximation of with $\At_{i\pm1/2}^n$ (\ref{DemiSomme}) is not appropriate to preserve the still
water steady state.
Thus, we have to find a suitable definition for $\At_{i\pm1/2}^n$ to  obtain an exactly well-balanced scheme. For this purpose, let us
assume \textbf{(SF)} and start with :

\noindent at the discrete level, the  still water steady state is perfectly maintained (see
\fig\ref{DiscretRepresentationDomain}): it exists $n$ such that for every $i$, if $Q_i^n=0$ and
$\forall i,$
\begin{description}
 \item[A1:] $\dsp c^2\ln\left(\dsp \frac{A_{i+1}^n}{S_{i+1}}\right)+g \mathcal{H}(\SE_{i+1}^n)\cos\theta +g b_{i+1} =
c^2\ln\left(\dsp \frac{A_{i}^n}{S_{i}}\right)+g \mathcal{H}(\SE_{i}^n)\cos\theta +g b_{i}$,
 \item[A2:] $\dsp AM_{i+1/2}^n = AP_{i-1/2}^n $,
 \item[A3:] $\dsp Q_{i+1/2}^n = Q_{i-1/2}^n $,
\end{description}
then, for all $l>n$ the conditions A1, A2 and A3 holds.

\noindent The condition A2 is satisfied if and only if
\begin{equation}\label{ConditionA2}
 {AM_{i+1/2}^n-AP_{i-1/2}^n} = \dsp \Delta A_{i+1/2}^n +
\frac{g}{2}\left(\frac{\At_{i+1/2}^n\psi_i^{i+1}}{\ctt_{i+1/2}^2}+\frac{\At_{i-1/2}^n\psi_{i-1}^{i}}{\ctt_{
i-1/2 }^2}\right) {=0}\,.
\end{equation}

\noindent The condition A3 is satisfied if and only if
\begin{equation}\label{ConditionA3}
\Delta Q_{i+1/2}^n =  \dsp \frac{g}{2}
\left\{
 \frac{\At_{i-1/2}^n\, \psi_{i-1}^{i}}{\ctt_{i-1/2}}-
 \frac{\At_{i+1/2}^n\, \psi_i^{i+1}}{\ctt_{i+1/2}}
\right\} + \frac{\Delta A_{i+1/2}^n}{2} \left(\ctt_{i-1/2} -
\ctt_{i+1/2}\right) = 0\,.
\end{equation}

\noindent The condition A1 is satisfied for pressurized flows if and only if
\begin{equation}\label{ConditionA1PF}
A_{i+1}^n = A_i^n \frac{S_{i+1}}{S_i}\exp\left(-\frac{g}{c^2}\left(\Delta b + \Delta R
\cos\theta \right)\right)\,.
\end{equation}

\begin{figure}[H]
\begin{center}
\includegraphics[scale=0.3]{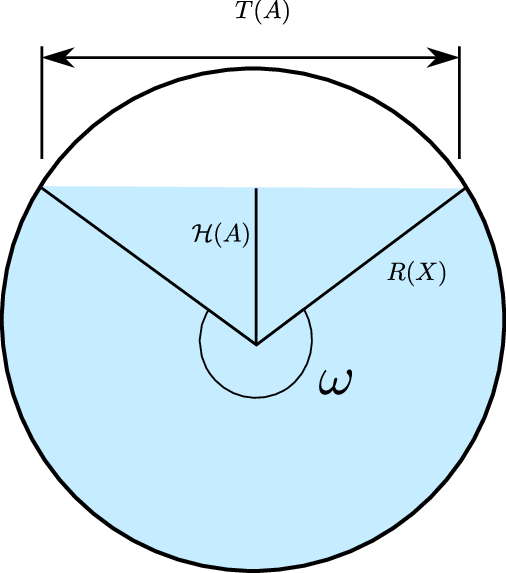}
\end{center}
\caption{angle $\omega$}
\label{ComputationOmega}
\end{figure}

\noindent The condition A1 is satisfied for free surface flows if and only if
\begin{equation}\label{ConditionA1FS}
h_{i+1}^n = \frac{h_i^n\cos\theta-\Delta b}{\cos\theta}\,.
\end{equation}
For circular cross-section pipe, $A_{i+1}^n$ is computed by:
\begin{equation}\label{ForCircularPipe}
A_{i+1}^n = \frac{R_{i+1}^2}{2} \left(\omega_{i+1}-\sin(\omega_{i+1})\right)
\end{equation}
with $\dsp \omega_{i+1} = \dsp   2\left(\pi -  \arccos\left(\frac{h_{i+1}}{R_{i+1}}\right)\right)  $ is the angle  displayed on
\fig\ref{ComputationOmega}.
\noindent Thus, the discrete  still
water steady state is perfectly maintained if and only if,  $(\At_{i-1/2}^n,\At_{i+1/2}^n)$ is the solution
of the non-linear system:
\begin{equation}\label{SNL}
\left\{
\begin{array}{lll}
 0&=& \dsp \Delta A_{i+1/2}^n
+\frac{g}{2}\left(\frac{\At_{i+1/2}^n\psi_i^{i+1}}{\ctt_{i+1/2}^2}+\frac{\At_{i-1/2}^n\psi_{i-1}^{i}}{\ctt_{i-1/2 }^2}\right)\\ & & \\
0&=& \dsp \frac{g}{2}\left\{ \frac{\At_{i-1/2}^n\, \psi_{i-1}^{i}}{\ctt_{i-1/2}}- \frac{\At_{i+1/2}^n\, \psi_i^{i+1}}{\ctt_{i+1/2}}\right\}
+ \frac{\Delta A_{i+1/2}^n}{2} \left(\ctt_{i-1/2} -
\ctt_{i+1/2}\right)
\end{array}
\right.
\end{equation}
\noindent where we have  replaced the expression of $A_{i+1}^n$ in (\ref{ConditionA2}-\ref{ConditionA3}) by (\ref{ConditionA1PF}) for
pressurized
by (\ref{ConditionA1FS}) for free surface flows:
\begin{equation*}
\Delta A_{i+1/2}^n =
\left\{
\begin{array}{lll}
\dsp A_i^n \left(\frac{S_{i+1}}{S_i}\exp\left(-\frac{g}{c^2}\left(\Delta b + \Delta R
\cos\theta \right)\right)-1 \right) &\textrm{ if } & E_i = 0\\
\dsp \mathcal{F}\left(\frac{h_i^n\cos\theta-\Delta b}{\cos\theta}\right)  &\textrm{ if } & E_i = 0
\end{array}
\right.\,.
\end{equation*}
with $\mathcal{F}: h\mapsto \mathcal{F}(h)=A$. For circular pipe, $\mathcal{F}$ is given by (\ref{ForCircularPipe}).

\noindent Finally, the numerical scheme (\ref{NumericalSchemeA}-\ref{NumericalSchemeQ})-(\ref{DefinitionOfTildebThetaSQ}) with $\At_{i\pm
1/2}^n$ as the solution of the non linear system (\ref{SNL}) defines an exactly well-balanced scheme.

\noindent For uniform pipe and pressurized flow, the previous system simply writes:
$$
\left\{
\begin{array}{lll}
 \dsp \Delta A_{i+1/2}^n +
\frac{g\,\Delta Z}{2\,c^2}\left(\At_{i+1/2}^n+\At_{i-1/2}^n\right) &=&0\\& &\\
\At_{i+1/2}^n&=&\At_{i-1/2}^n
\end{array}
\right.\,.
$$
The solution is easily obtained by:
\begin{equation}\label{AtildeByAi}
\At_{i+1/2}^n = -\frac{c^2}{g}\frac{\Delta A_{i+1/2}^n}{\Delta b} =  \dsp
-\frac{c^2}{g}\frac{A_i^n\left(\exp\left(\dsp-\frac{g}{c^2}\Delta b \right)-1\right)}{\Delta
b}\,.
\end{equation}
Let us also remark that  using the relation $A_i^n =
A_{i+1}^n\exp\left(\dsp\frac{g}{c^2}\Delta b \right)$, we have:
\begin{equation}\label{AtildeByAi+1}
\At_{i+1/2}^n = -\frac{c^2}{g}\frac{\Delta A_{i+1/2}^n}{\Delta b} =  \dsp
-\frac{c^2}{g}\frac{A_{i+1}^n\left(1-\exp\left(\dsp\frac{g}{c^2}\Delta b \right)\right)}{\Delta
b}\, .\end{equation}
It follows that $\At_{i+1/2}^n$ can be expressed as the mean value of  (\ref{AtildeByAi}) and (\ref{AtildeByAi+1}) as follows:
$$\At_{i+1/2}^n = \dsp-\frac{c^2}{g\,\Delta
b}\left\{\frac{ A_{i+1}^n\left(1-\exp\left(\dsp\frac{g}{c^2}\Delta b \right)\right)
+ A_i^n\left(\exp\left(\dsp-\frac{g}{c^2}\Delta b \right)-1\right)}{2}\right\}\,.$$
For small $\Delta X$, we have $$\At_{i+1/2}^n\approx \frac{A_i^n+A_{i+1}^n}{2}\,.$$ It follows that the scheme
(\ref{NumericalSchemeA}-\ref{NumericalSchemeQ})-(\ref{DefinitionOfTildebThetaSQ})
with (\ref{DemiSomme}) is the zero order approximation of the solution given by the EWBS.

\noindent The same analysis shows that the free surface numerical scheme with (\ref{DemiSomme}) is also the zero order approximation of
the solution given by the EWBS.

\begin{figure}[H]
\begin{center}
\includegraphics[scale=1.0]{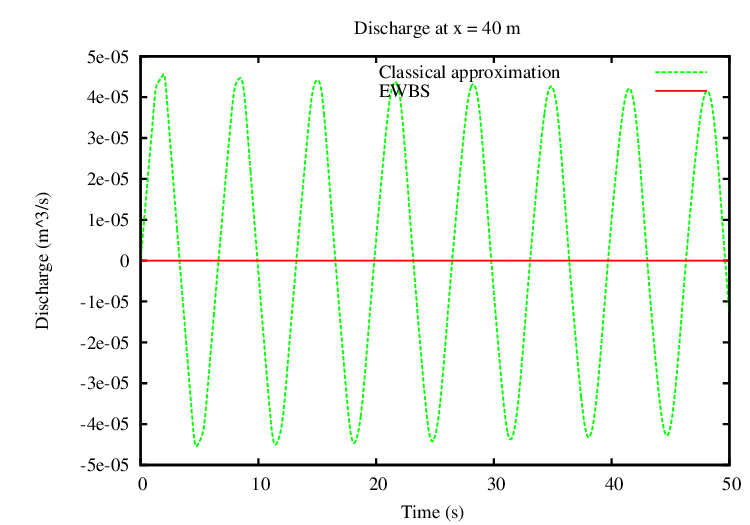}

\includegraphics[scale=1.0]{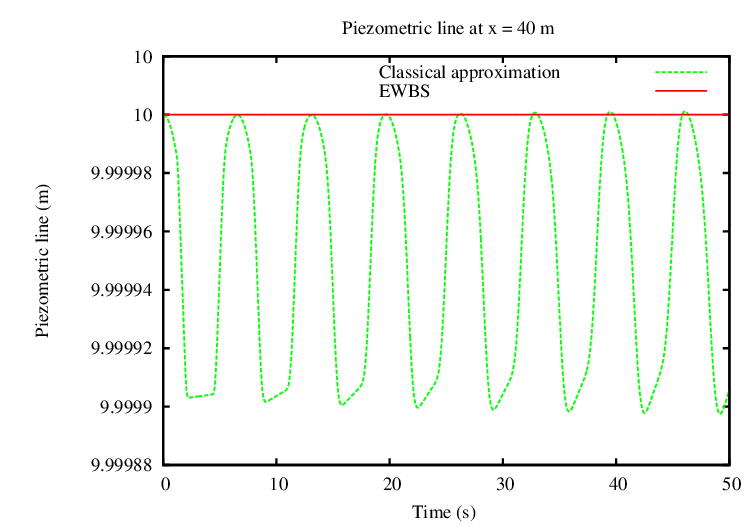}
\end{center}
\caption{The numerical scheme (\ref{NumericalSchemeA}-\ref{NumericalSchemeQ})-(\ref{DefinitionOfTildebThetaSQ}) with (\ref{DemiSomme}) and
the EWBS for pressurized   still
water steady state with  $c=30$.}
\label{ComparaisonMeanMatrixStationnaryMatrixc30}
\end{figure}

\begin{figure}[H]
\begin{center}
\includegraphics[scale=1.0]{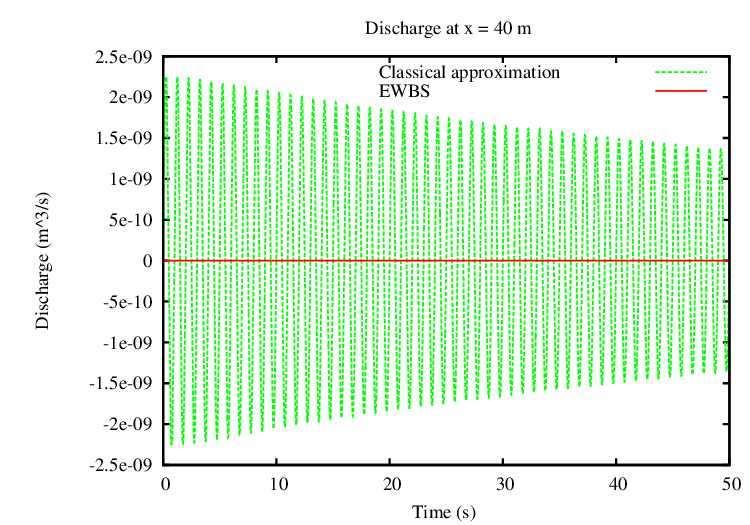}

\includegraphics[scale=1.0]{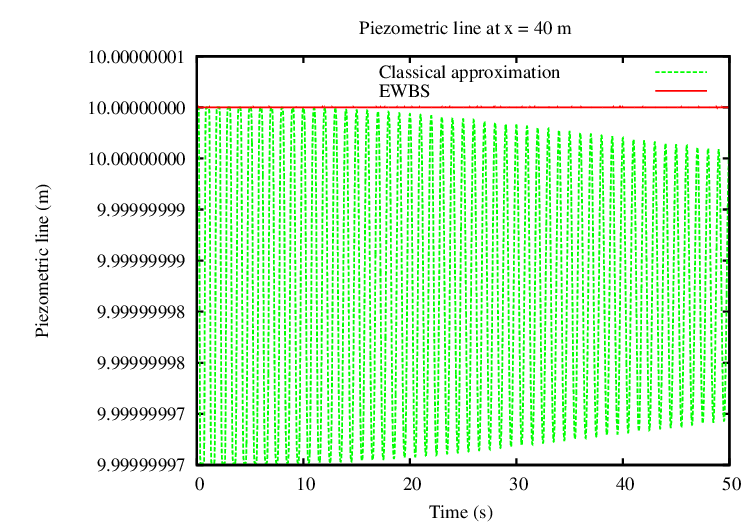}
\end{center}
\caption{The numerical scheme (\ref{NumericalSchemeA}-\ref{NumericalSchemeQ})-(\ref{DefinitionOfTildebThetaSQ}) with (\ref{DemiSomme}) and
the EWBS for pressurized  still
water steady state with    $c=200$.}
\label{ComparaisonMeanMatrixStationnaryMatrixc100}
\end{figure}

\begin{figure}[H]
\begin{center}
\includegraphics[scale=1.0]{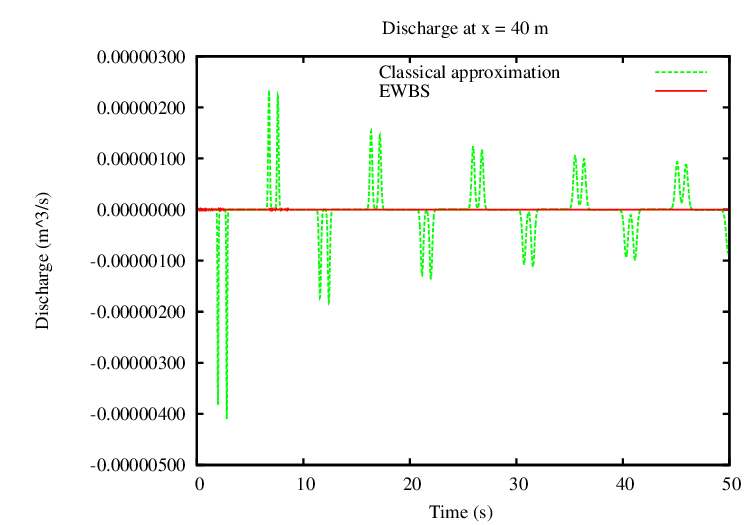}

\includegraphics[scale=1.0]{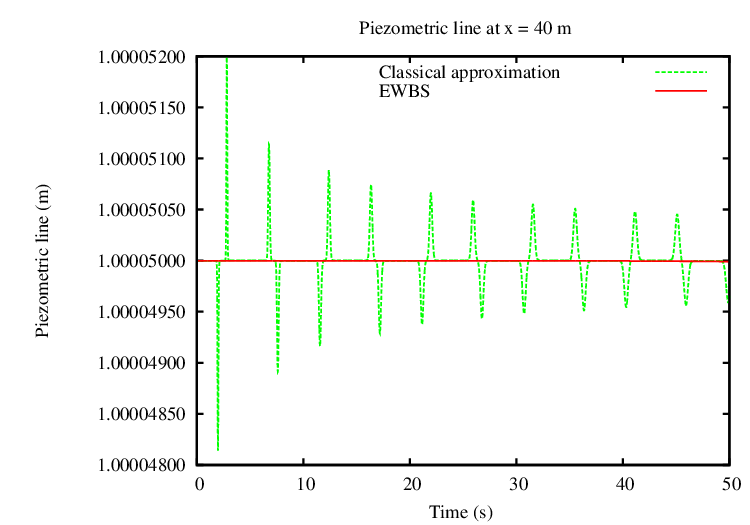}
\end{center}
\caption{The numerical scheme (\ref{NumericalSchemeA}-\ref{NumericalSchemeQ})-(\ref{DefinitionOfTildebThetaSQ}) with (\ref{DemiSomme}) and
the EWBS for free surface  still
water steady state.}
\label{StatioSL}
\end{figure}

\paragraph*{}
On \fig\ref{ComparaisonMeanMatrixStationnaryMatrixc30}, \fig\ref{ComparaisonMeanMatrixStationnaryMatrixc100}, \fig\ref{StatioSL}, we
display the  still water steady pressurized and free surface state computed by the EWBS and the scheme with the approximation
(\ref{DemiSomme}). The steady state with the EWBS is numerically well preserved while as pointed out before (see Subsection
\ref{SubSectionSteadyStatesWithClassicalApproximation}) the classical approach is not convenient.

\begin{figure}[H]
\begin{center}
\includegraphics[scale=1.0]{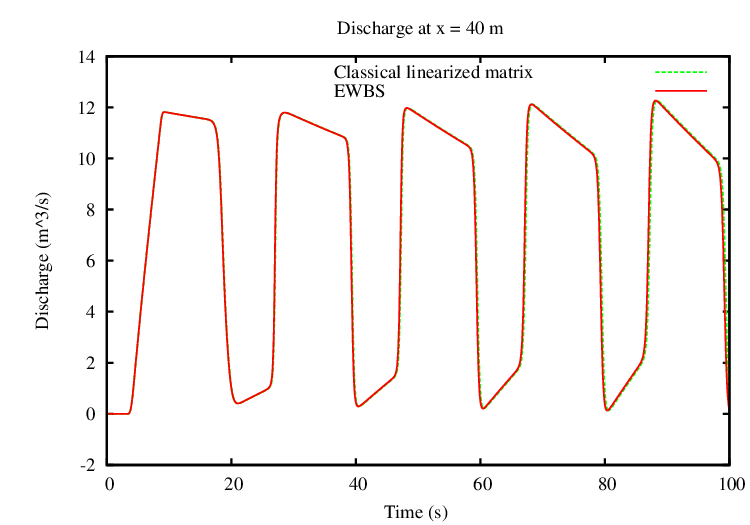}

\includegraphics[scale=1.0]{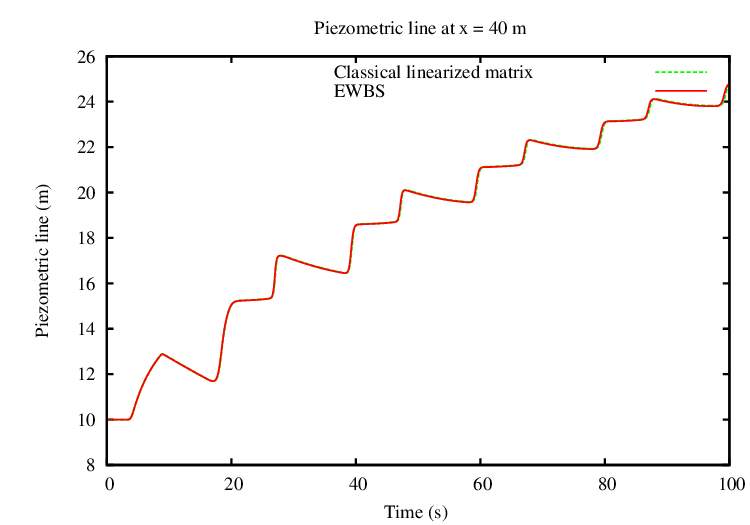}
\end{center}
\caption{A non stationary test to compare the EWBS and the numerical scheme
(\ref{NumericalSchemeA}-\ref{NumericalSchemeQ})-(\ref{DefinitionOfTildebThetaSQ}) with (\ref{DemiSomme}).}
\label{TestNonStationnaire}
\end{figure}

\noindent We also display an unsteady simulation on \fig\ref{TestNonStationnaire} where the results of the two methods are very well
reproduced.

\subsection{Remarks concerning  mixed  still water steady state}
The previous sections deals with the well-balanced property of the numerical scheme
(\ref{NumericalSchemeA}-\ref{NumericalSchemeQ})-(\ref{DefinitionOfTildebThetaSQ}) for
free surface and pressurized flows. To use the well-balanced scheme developed in Subsections
(\ref{SubSectionSteadyStatesWithClassicalApproximation}) and
(\ref{SubSectionEWBS}), we
start from the discrete representation of a mixed  still
water steady state (as displayed on \fig\ref{DiscretRepresentationDomain}).
Assume that there  exists $i_0$ such that all cells $m_i$ on the left hand side of the interface $x_{i_0+1/2}$ are free
surface while the other are pressurized:
$$
\left\{
\begin{array}{lll}
E_i = 0  & \textrm{ if }& i\leqslant i_0 \\
E_i = 1  & \textrm{ if }& i> i_0
\end{array}\,.
\right.
$$

\begin{figure}[H]
\begin{center}
\includegraphics[scale=0.4]{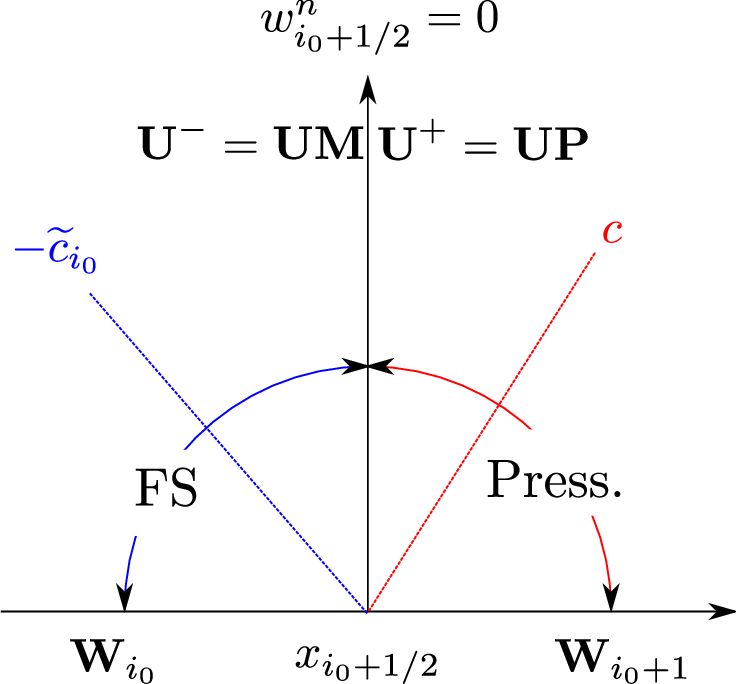}
\end{center}
\caption{Mixed free surface-pressurized  still
water steady state}
\label{Wzero}
\end{figure}

\noindent The interface $x_{i_0+1/2}$ is such that the speed of propagation of the interface is null:
$w_{i_0+1/2}^n=\dsp\frac{Q_{i+1}^n-Q_i^n}{A_{i+1}^n-A_i^n}=0$ (see \fig\ref{Wzero}).
Therefore, $\UM = \U^-$ and $\UP = \U^+$.
\noindent For example on \fig\ref{Wzero}, we apply the free surface numerical scheme on the left hand side of the interface  and the
pressurized one on the right cells . As the EWBS preserves  the pressurized and free surface still water steady state,  it also preserves
the mixed  still water steady state.

\section{Numerical tests}\label{NumericalTests}
%
The numerical validation for pipes with constant section and slope  has been previously studied by two of the authors in \cite{BG07,BG08}
and thus are not presented in this paper.
Since experimental data for mixed flows in any pipe are not available, we focus on the behavior of our method for several circular
cross-section contracting and expanding pipe. Notice that, the equivalent pipe method is not relevant for the mixed flows as pointed out by
\cite{A03,SWB98,WS78} for instance.

\noindent The mixed flow case is numerically performed on a water hammer test. Starting from an horizontal free surface  still
water steady state, the
water hammer occurs immediately after the increase of the upstream piezometric head while the downstream  discharge is set to $0$. The
prescribed hydrograph produces a travelling wave which produces a pressurized state propagating from upstream to downstream end. Physically
an trapped air pocket may appear: it is not taken into account in the \textbf{PFS}-model. Actually, the trapped air pockets vanish or move; some
parts of these pockets undergo condensation/vaporisation and others parts move and lead to a two phase flow. Consequently the sound speed
decreases. As our model does not take into account  these phenomena, the value of $c$ is assumed to be constant. Moreover we should have to
deal with the entrapment of air bubbles which have a non negligible effect (see \cite{HM82,S76} for instance).

\begin{figure}[H]
\centering
\includegraphics[angle=0,scale=1.0]{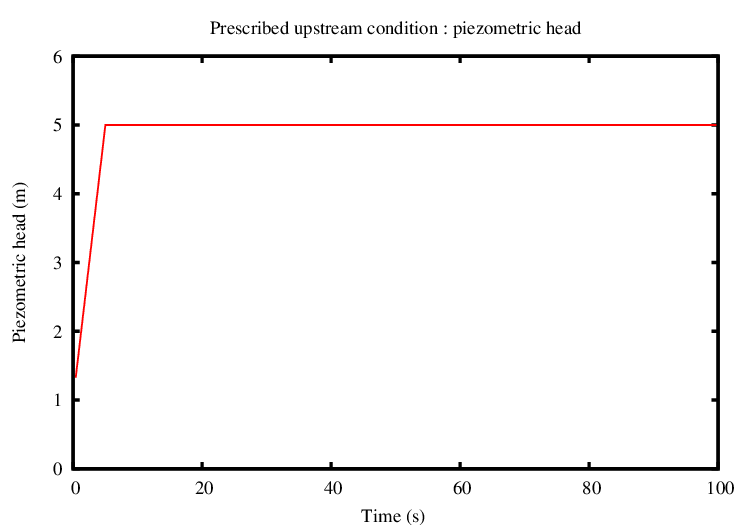}

\includegraphics[angle=0,scale=1.0]{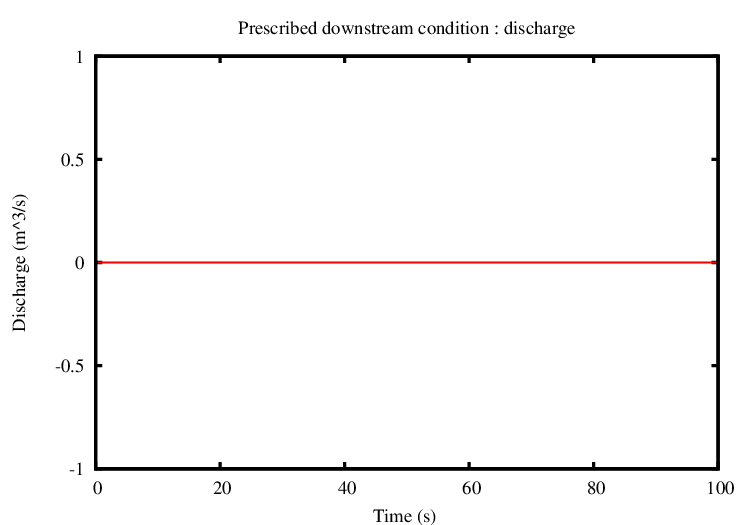}
\caption{Boundary conditions.}\label{BoundaryCondition}
\end{figure}

\noindent The numerical experiments are performed in the case of a $100$ $m$ long closed circular pipe    at altitude  $b_0 = 1\,m$
with $0$ slope which corresponds to the elevation and slope of the main pipe axis (we have $Z = b(X) = 0,\,\forall X$). The Manning
roughness coefficient is $1/K_s= 0.012 \; s/m^{1/3}$. The simulation starts from a
steady state as a free surface flow with a discharge $Q = 0 \;m^3/s$.
The upstream boundary condition  is a prescribed hydrograph (see \fig\ref{BoundaryCondition}) while the downstream discharge is kept
constant to $0\,m^3/s$ (as displayed on \fig\ref{BoundaryCondition}).
We compare then the results obtained for uniform, contracting an expanding pipes. For each test, the parameters are the same except the
downstream diameter: the upstream diameter is kept constant to $D=1\,m$. The contracting pipe is chosen  for  $D =0.6\,m$ and the
expanding one for  $D = 1.4\,m$ (where $D$ denotes the downstream diameter).
\noindent Let us recall that the zero water level corresponds to the main pipe axis. The piezometric head is defined  as $z + p$:
$$\dsp \left\{\begin{array}{l}
\dsp p = 2R+\frac{c^2 \, (\rho -\rho_0)}{\rho_0 \, g}
\mbox{ if the flow is pressurized}\\
p = h \mbox{ the water height if the flow is free surface}
\end{array}
\right.
$$

\begin{figure}[H]
\begin{center}
\includegraphics[scale=1.0]{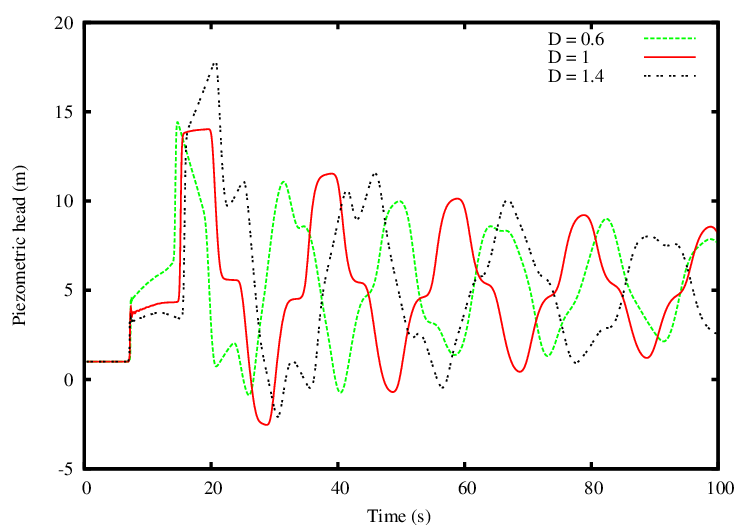}

\includegraphics[scale=1.0]{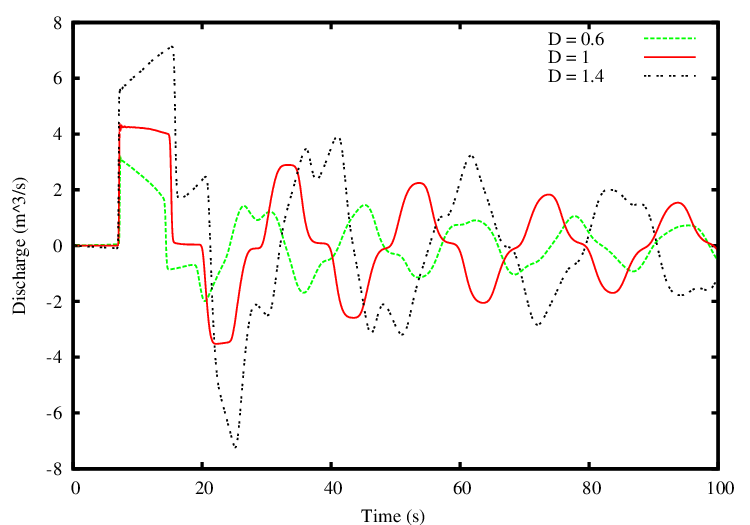}
\end{center}
\caption{Piezometric head and discharge at $X=50m$.}
\label{X50}
\end{figure}

\begin{figure}[H]
\centering
\includegraphics[angle=0,scale=0.75]{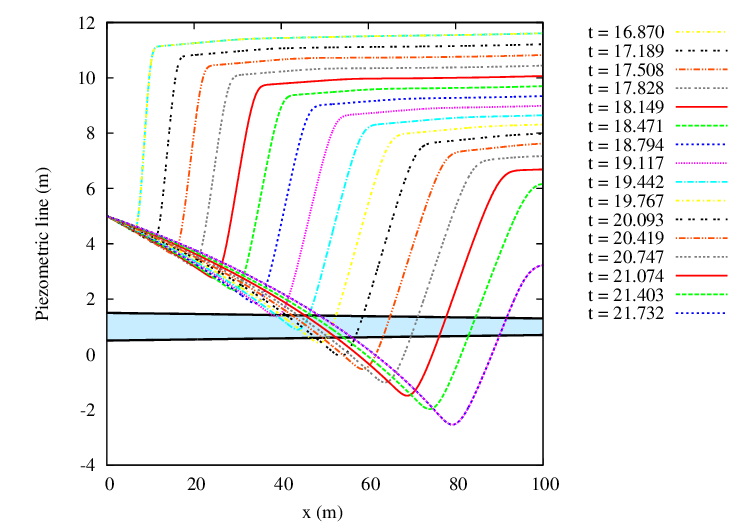}

\includegraphics[angle=0,scale=0.75]{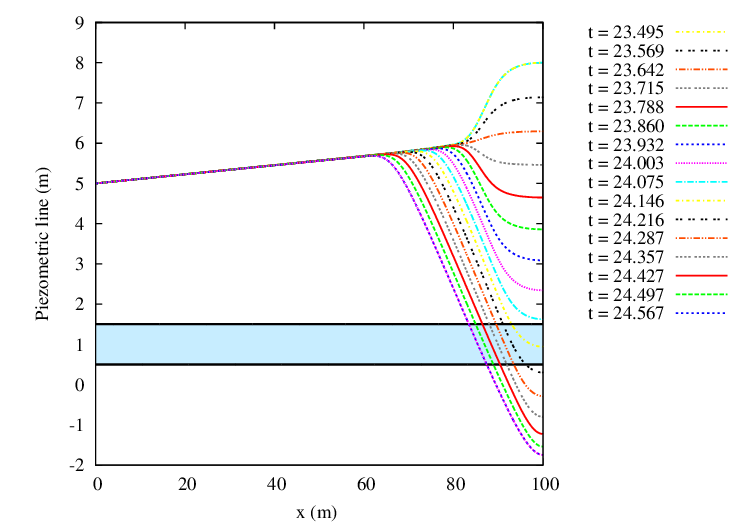}

\includegraphics[angle=0,scale=0.75]{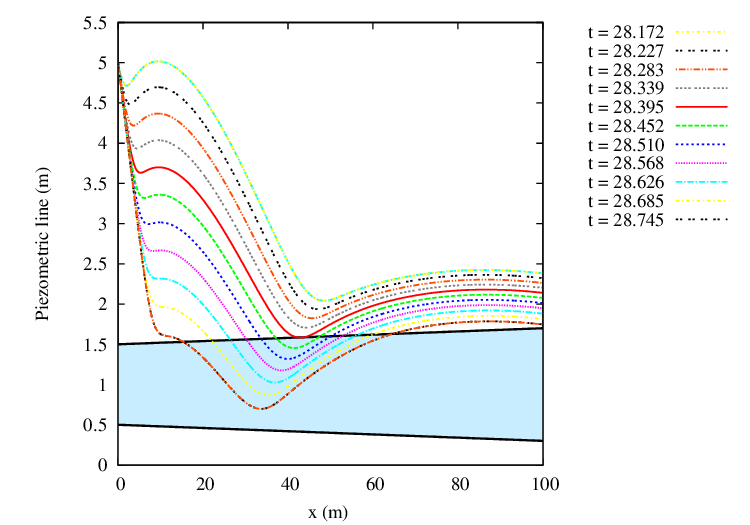}

\caption{ Observation of the depression localised approximatively at time $t = 19.117$ (contracting pipe), $t = 24.075$ (uniform pipe) and $t = 28.395$ (expanding pipe).  }
\label{Depress}
\end{figure}

\begin{figure}[H]
\centering
\includegraphics[angle=0,scale=0.75]{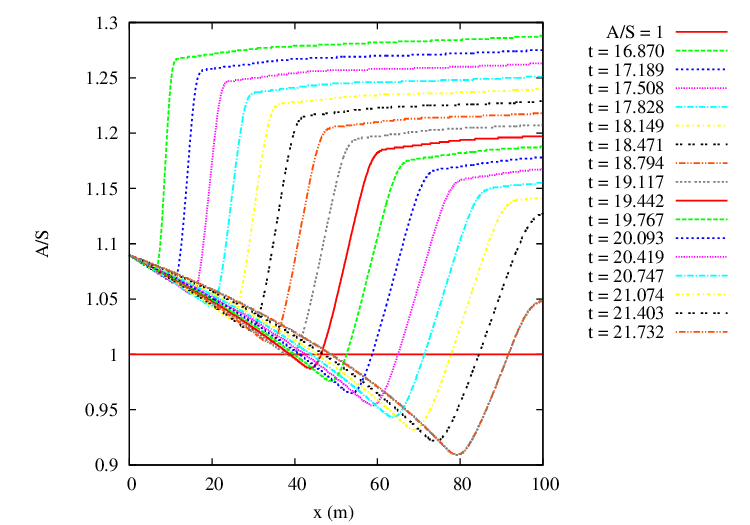}

\includegraphics[angle=0,scale=0.75]{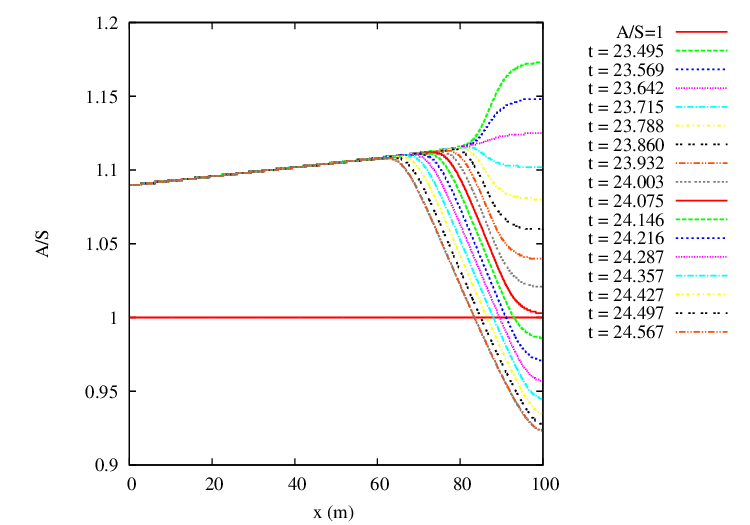}

\includegraphics[angle=0,scale=0.75]{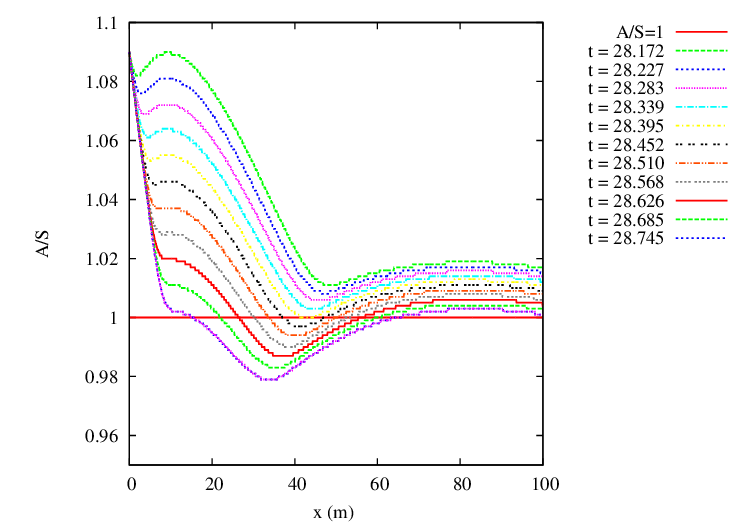}

\caption{ Observation of the depression localised approximatively at time $t = 19.117$ (contracting pipe), $t = 24.075$ (uniform pipe) and $t = 28.395$ (expanding pipe).  }
\label{DepressX}
\end{figure}

\noindent Results are then represented on \fig\ref{X50}. The sudden raise of the upstream piezometric level produces a  pressurized
state with a travelling wave.  A  water hammer is then observed since the downstream discharge is null.
A careful analysis of the flow which is performed by the variable $E$ or equivalently by $\dsp\frac{A}{\SE}$ (see Remark
\ref{RemarkDepressionOverpressureStates} and \fig{\ref{RhoCriterion}})   shows that after this transition point, the flow is pressurized but
in depression which starts approximatively at time $19s$ for the contracting pipe, $24s$ for the uniform pipe and $28s$ for the expanding
one. We display the piezometric line for different times around the depression time (see \fig\ref{Depress}) and also the graph of the
function $\dsp\frac{A}{\SE}$ which  confirm that the observed times correspond exactly to a depression state for each pipes (see
\fig\ref{DepressX}).

\noindent We also observe  a little smoothing effect and absorption due to the
first order discretisation type.

\section{Conclusion}
We have  derived a free surface and a pressurized  model which have been coupled using a common set of variables  and a suitable pressure
law. We have thus obtained a mathematical model for unsteady mixed flows in non uniform water pipes, that we have called \textbf{PFS}-model.
This  model takes into account the local perturbation of the section and of the slope.

The \textbf{PFS} model is numerically solved by a  VFRoe scheme  using the interfacial upwind to include the source terms into  the
numerical fluxes. We have shown that the classical approximation of the convection matrix (the Roe matrix approximation) is not
suitable to preserve the  still
water steady state (except for the pressurized case where the value of $c$ helps the scheme  to maintain this  state). Moreover, we
have
proposed a manner to obtain an exactly well-balanced scheme.

As mentioned in \cite{BG07} this numerical method with the classical approximation of the convection matrix, for constant section,
reproduces correctly  laboratory tests for uniform
pipes and can deal with multiple points of transition between the two types of flows. As pointed out before, due to the lack of
experimental data for non uniform pipes, we have only
shown the behavior of the piezometric line which seems reasonable (at less no major difference was observed).

As a well-known feature on approximate Godunov scheme, the upwinding of the source terms may introduce stationary waves with  a vanishing
denominator when critical flows occurs i.e. $u \approx c$. Moreover, in its actual form, the presented numerical scheme is not able to deal with drying and
flooding area. Nevertheless, it may be possible to introduce a cut-off function to avoid division by zero for each problems: critical
stationary waves, drying area and flooding area. But it is not the better choice that we can propose, since,  truncation of
the wet
area induces a loss of water mass leading to the non-conservativity of the mass. Nevertheless, at the present time, we are
interested in a mathematical kinetic formulation of the \textbf{PFS} model and the construction of a numerical
kinetic scheme that avoids all these drawbacks \cite{BEG09_4}.

The next step is to take into account the air entrainment which may have non negligible effects on the behaviour of the
piezometric head. A first approach has been derived in the case of perfect fluid and perfect gas modelised as a bilayer model based on the
\textbf{PFS}-model \cite{BEG09_3}.

\vspace*{0.5cm}

\noindent \textbf{Aknowlegments:} The authors  wish to thank the two referees for their  careful reading  of the first version of the article and useful remarks.
\bibliographystyle{plain}


\newpage

\begin{figure}[H]
\centering
\includegraphics[angle=0,scale=0.75]{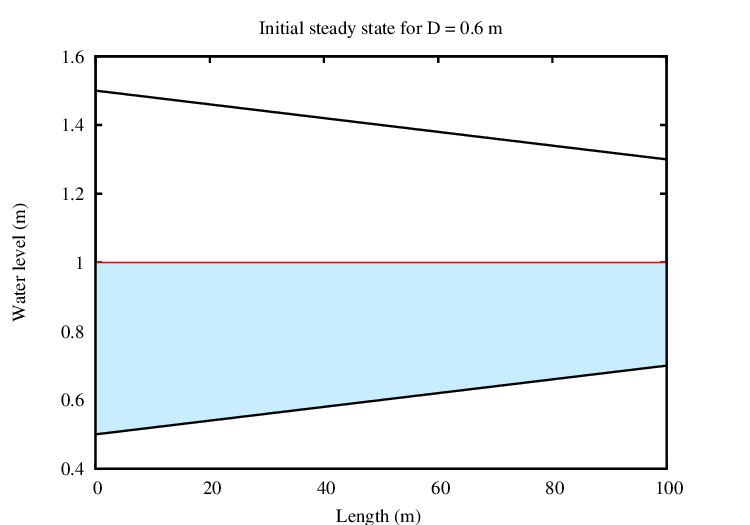}

\includegraphics[angle=0,scale=0.75]{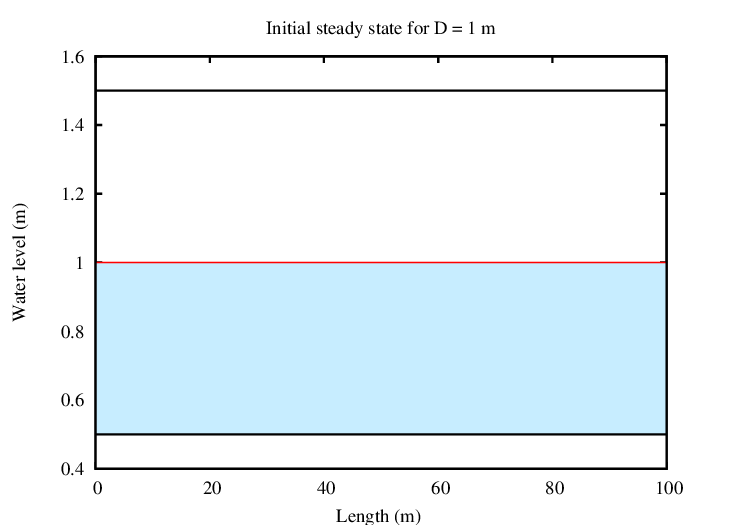}

\includegraphics[angle=0,scale=0.75]{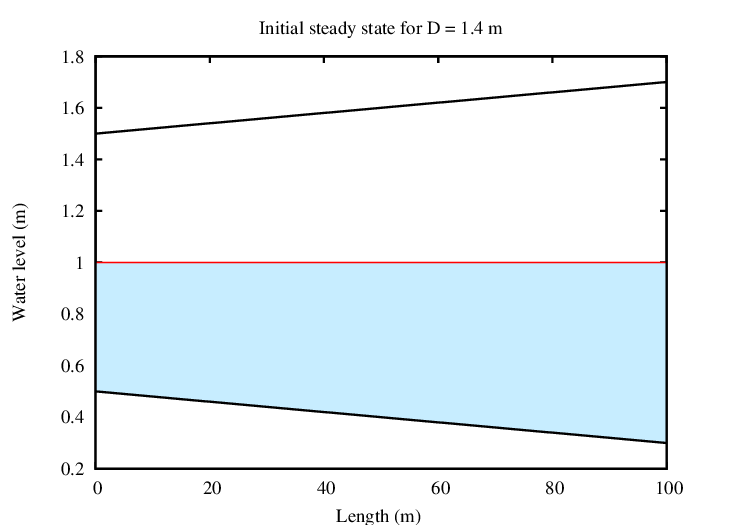}
\caption{Initial still water steady state for contracting, uniform and expanding pipes  \label{statio}.}
\end{figure}

\end{document}